\documentclass[12pt,a4paper,twoside]{article}

\newcommand{\pageformat}[6]{\setlength{\hoffset}{-1in}
                  \setlength{\voffset}{-1in}
                  \addtolength{\hoffset}{#5}
                            \addtolength{\voffset}{#6}
                            \setlength{\oddsidemargin}{#1}
                            \setlength{\evensidemargin}{#2}
                            \setlength{\textwidth}{\paperwidth}
                  \addtolength{\textwidth}{-\oddsidemargin}
                  \addtolength{\textwidth}{-\evensidemargin}
                  \addtolength{\textwidth}{-\marginparsep}
                  \addtolength{\textwidth}{-\marginparwidth}
                            \setlength{\topmargin}{#3}
                            \setlength{\textheight}{\paperheight}
                  \addtolength{\textheight}{-\topmargin}
                  \addtolength{\textheight}{-\headheight}
                  \addtolength{\textheight}{-\headsep}
                  \addtolength{\textheight}{-\footskip}
                  \addtolength{\textheight}{-#4}}
\pageformat{2cm}{3cm}{24mm}{24mm}{1pt}{0pt}

\usepackage{ifthen}
\newboolean{article}
    \setboolean{article}{true}
\newboolean{report}
\newboolean{book}
\newboolean{letter}
\newboolean{german}
\newboolean{italian}
\newboolean{nobaselinestretch}
\newboolean{nosectionappendix}
\newboolean{oldtoc}
\newboolean{nosectionequn}
\newboolean{notheorem}

\ifthenelse{\boolean{german}}{
    \usepackage{german}}{}

\usepackage[latin1]{inputenc}

\usepackage{amsmath}
\usepackage{amssymb}
\usepackage[mathscr]{eucal}

\ifthenelse{\boolean{notheorem}}{}{
    \usepackage{theorem}}

\ifthenelse{\boolean{nobaselinestretch}}{}{
    \renewcommand{\baselinestretch}{1.25}}

\newenvironment{env}[2]{\begin{#1}#2\end{#1}}{}
    \newcommand{\beq}[1]{\begin{env}{equation}{#1}}
    \newcommand{\beqn}[1]{\begin{env}{equation*}{#1}}
    \newcommand{\bal}[1]{\begin{env}{align}{#1}}
    \newcommand{\baln}[1]{\begin{env}{align*}{#1}}
    \newcommand{\bga}[1]{\begin{env}{gather}{#1}}
    \newcommand{\bgan}[1]{\begin{env}{gather*}{#1}}
    \newcommand{\bflal}[1]{\begin{env}{flalign}{#1}}
    \newcommand{\bflaln}[1]{\begin{env}{flalign*}{#1}}
    \newcommand{\bmu}[1]{\begin{env}{multline}{#1}}
    \newcommand{\bmun}[1]{\begin{env}{multline*}{#1}}
    \newcommand{\bsp}[1]{\begin{env}{split}{#1}}

    \newcommand{\eeq}{\end{env}}
    \newcommand{\eeqn}{\end{env}}
    \newcommand{\eal}{\end{env}}
    \newcommand{\ealn}{\end{env}}
    \newcommand{\ega}{\end{env}}
    \newcommand{\egan}{\end{env}}
    \newcommand{\eflal}{\end{env}}
    \newcommand{\eflaln}{\end{env}}
    \newcommand{\emu}{\end{env}}
    \newcommand{\emun}{\end{env}}
    \newcommand{\esp}{\end{env}}

\newcommand{\lf}{\vspace{2ex}}

\renewcommand{\bf}[1]{\textbf{#1}}
\renewcommand{\it}[1]{\textit{#1}}

\renewcommand{\sf}[1]{\textsf{#1}}

\renewcommand{\tt}[1]{\texttt{#1}}
\newcommand{\hl}[1]{\bf{\it{#1}}}

\newcommand{\mbf}[1]{\mathbf{#1}}
\newcommand{\msf}[1]{\text{\small$\sf{#1}$}}

\newcommand{\cmc}[1]{\mathcal{#1}}
\newcommand{\eus}[1]{\mathscr{#1}}
\newcommand{\euf}[1]{\mathfrak{#1}}
\newcommand{\bb}[1]{\mathbb{#1}}

\newcommand{\nbd}[1]{$#1$\nobreakdash--}
\newcommand{\ol}[1]{\overline{#1}}

\newcommand{\wt}[1]{\widetilde{#1}}
\newcommand{\wh}[1]{\widehat{#1}}
\newcommand{\ve}{\varepsilon}
\newcommand{\vt}{\vartheta}
\newcommand{\vk}{\varkappa}
\newcommand{\vp}{\varphi}

\newcommand{\Om}{\Omega}
\newcommand{\abs}[1]{\left\lvert#1\right\rvert}
\newcommand{\norm}[1]{\left\lVert#1\right\rVert}

\newcommand{\bnorm}[1]{\bigl\lVert#1\bigr\rVert}

\newcommand{\snorm}[1]{\norm{\smash{#1}}}

\newcommand{\bfam}[1]{\bigl(#1\bigr)}
\newcommand{\Bfam}[1]{\Bigl(#1\Bigr)}
\newcommand{\AB}[1]{\langle#1\rangle}

\newcommand{\CB}[1]{\{#1\}}
\newcommand{\bCB}[1]{\bigl\{#1\bigr\}}
\newcommand{\BCB}[1]{\Bigl\{#1\Bigr\}}
\newcommand{\SB}[1]{[#1]}

\newcommand{\LO}[1]{(#1]}

\newcommand{\set}[2][]{
    \ifthenelse{\equal{#1}{}}{
        \CB{#2}}{
        \CB{#1~|~#2}}}
\newcommand{\bset}[2][]{
    \ifthenelse{\equal{#1}{}}{
        \bCB{#2}}{
        \bCB{#1~|~#2}}}
\newcommand{\Bset}[2][]{
    \ifthenelse{\equal{#1}{}}{
        \BCB{#2}}{
        \BCB{#1~\big|~#2}}}

\DeclareMathOperator{\ls}{\normalfont\msf{span}}
\DeclareMathOperator{\cls}{\ol{\ls}}

\DeclareMathOperator{\id}{\normalfont\msf{id}}

\renewcommand{\ker}{\operatorname{\msf{ker}}}

\newcommand{\C}{\bb{C}}

\newcommand{\N}{\bb{N}}

\newcommand{\R}{\bb{R}}

\newcommand{\cA}{\cmc{A}}
\newcommand{\cB}{\cmc{B}}
\newcommand{\cC}{\cmc{C}}

\newcommand{\cT}{\cmc{T}}

\newcommand{\sL}{\eus{L}}

\newcommand{\sT}{\eus{T}}

\newcommand{\eR}{\euf{R}}
\newcommand{\eS}{\euf{S}}

\newcommand{\eZ}{\euf{Z}}
\newcommand{\U}{\mbf{1}}

\newcommand{\G}{\Gamma}

\newcommand{\I}{{I\!\!\!\;I}}

\newcommand{\RL}{{\rm L}}

\ifthenelse{\boolean{nosectionequn}}{}{
    \numberwithin{equation}{section}
    }

\ifthenelse{\boolean{article}\or\boolean{letter}\or\boolean{nosectionequn}}{
    \setboolean{nosectionappendix}{true}}{}
\ifthenelse{\boolean{nosectionappendix}}{}{
    \renewcommand{\appendix}{
        \chapter*{\appendixname}
        \addcontentsline{toc}{chapter}{\appendixname}
        \renewcommand{\thesection}{\Alph{section}}
        \setcounter{section}{0}}}
   
\ifthenelse{\boolean{report}\or\boolean{book}}{
    }{}

\ifthenelse{\boolean{notheorem}}{}{
        \newcommand{\mnname}{Mathematical note.}
        \newcommand{\enname}{End of the note.}
        \newcommand{\definame}{Definition.}
        \newcommand{\propname}{Proposition.}
        \newcommand{\lemname}{Lemma.}
        \newcommand{\exname}{Example.}
        \newcommand{\exername}{Exercise.}
        \newcommand{\remname}{Remark.}
        \newcommand{\obname}{Observation.}
        \newcommand{\thmname}{Theorem.}
        \newcommand{\corname}{Corollary.}
        \newcommand{\proofname}{Proof.}
    \ifthenelse{\boolean{german}}{
        \renewcommand{\mnname}{Mathematische Notiz.}
        \renewcommand{\enname}{Ende der Notiz.}
        \renewcommand{\exname}{Beispiel.}
        \renewcommand{\exername}{Übung.}
        \renewcommand{\remname}{Bemerkung.}
        \renewcommand{\obname}{Beobachtung.}
        \renewcommand{\thmname}{Satz.}
        \renewcommand{\corname}{Korollar.}
        \renewcommand{\proofname}{Beweis.}}{}
    \ifthenelse{\boolean{italian}}{
        \renewcommand{\mnname}{Nota matematica.}
        \renewcommand{\enname}{Fina della nota.}
        \renewcommand{\definame}{Definizione.}
        \renewcommand{\propname}{Proposizione.}
        \renewcommand{\exname}{Esempio.}
        \renewcommand{\exername}{Esercizio.}
        \renewcommand{\remname}{Nota.}
        \renewcommand{\obname}{Osservazione.}
        \renewcommand{\thmname}{Teorema.}
        \renewcommand{\corname}{Corollario.}
        \renewcommand{\proofname}{Dimostrazione.}

       \renewcommand{\appendixname}{Appendice}

       }{}
    \theoremheaderfont{\normalfont\bfseries}
    \theoremstyle{change}
        \theorembodyfont{\rmfamily}
            \newtheorem{emp}{}[section]
                \newcommand{\bemp}[1][]{
                    \begin{emp}\hskip-\labelsep\bf{#1}\hskip\labelsep}
                \newcommand{\eemp}{\end{emp}}
\newtheorem{itemp}[emp]{}
                \newcommand{\bitemp}[1][]{
                    \begin{itemp}\hskip-\labelsep\bf{#1}\hskip\labelsep\normalfont\itshape}
                \newcommand{\eitemp}{\end{itemp}}
            \newtheorem{mn}[emp]{\mnname}
                \newcommand{\bnm}{\begin{mn}~\begin{quotation}\renewcommand{\baselinestretch}{1}\small\noindent\ignorespaces}
                \newcommand{\enm}{\end{quotation}\hfill\bf{\enname}\end{mn}}
            \newtheorem{ex}[emp]{\exname}
                \newcommand{\bex}{\begin{ex}}
                \newcommand{\eex}{\end{ex}}
            \newtheorem{exer}[emp]{\exername}
                \newcommand{\bexer}{\begin{exer}}
                \newcommand{\eexer}{\end{exer}}
            \newtheorem{defi}[emp]{\definame}
                \newcommand{\bdefi}{\begin{defi}}
                \newcommand{\edefi}{\end{defi}}
            \newtheorem{rem}[emp]{\remname}
                \newcommand{\brem}{\begin{rem}}
                \newcommand{\erem}{\end{rem}}
            \newtheorem{ob}[emp]{\obname}
                \newcommand{\bob}{\begin{ob}}
                \newcommand{\eob}{\end{ob}}
        \theorembodyfont{\normalfont\itshape}
            \newtheorem{thm}[emp]{\thmname}
                \newcommand{\bthm}{\begin{thm}}
                \newcommand{\ethm}{\end{thm}}
            \newtheorem{prop}[emp]{\propname}
                \newcommand{\bprop}{\begin{prop}}
                \newcommand{\eprop}{\end{prop}}
            \newtheorem{cor}[emp]{\corname}
                \newcommand{\bcor}{\begin{cor}}
                \newcommand{\ecor}{\end{cor}}
            \newtheorem{lem}[emp]{\lemname}
                \newcommand{\blem}{\begin{lem}}
                \newcommand{\elem}{\end{lem}}
\newenvironment{empn}[1]{\lf\noindent\bf{#1}\ignorespaces\hskip\labelsep}{\lf}
		\newcommand{\bempn}[1]{\begin{empn}{#1}}
		\newcommand{\eempn}{\end{empn}}
		\newcommand{\bitempn}[1]{\begin{empn}{#1}\normalfont\itshape}
		\newcommand{\eitempn}{\end{empn}}
                \newcommand{\bnmn}{\begin{empn}{\mnname}~\begin{quotation}\renewcommand{\baselinestretch}{1}\small\noindent\ignorespaces}
                \newcommand{\enmn}{\end{quotation}\hfill\bf{\enname}\end{empn}}
		\newcommand{\bexn}{\begin{empn}{\exname}}
		\newcommand{\eexn}{\end{empn}}
		\newcommand{\bexern}{\begin{empn}{\exername}}
		\newcommand{\eexern}{\end{empn}}
		\newcommand{\bdefin}{\begin{empn}{\definame}}
		\newcommand{\edefin}{\end{empn}}
		\newcommand{\bremn}{\begin{empn}{\remname}}
		\newcommand{\eremn}{\end{empn}}
		\newcommand{\bobn}{\begin{empn}{\obname}}
		\newcommand{\eobn}{\end{empn}}

\newcommand{\qedsymbol}{~\rule[-0.35mm]{2mm}{2mm}}
    \newcounter{proof}[emp]
    \newenvironment{Proof}[1]{
        \vspace{1ex}
        \renewcommand{\item}[1][\stepcounter{proof}(\roman{proof})]%
            {##1\hskip\labelsep}
        \noindent\textsc{#1\hskip\labelsep}}{
        \nolinebreak\qedsymbol}
    \newcommand{\proof}[1][\proofname]{
        \begin{Proof}{#1}\ignorespaces}
    \newcommand{\qed}{\end{Proof}}
    \newcommand{\noqed}{
        \renewcommand{\qedsymbol}{}
        \end{Proof}}}
    \ifthenelse{\boolean{italian}}{
        \renewcommand{\proofname}{Dimostrazione.}}{}

\newcommand{\GL}{\Lambda}
\newcommand{\ga}{\alpha}
\newcommand{\GG}{\Gamma}
\newcommand{\GD}{\Delta}
\newcommand{\SA}{\cA}
\newcommand{\ot}{\otimes}
\newcommand{\gd}{\delta}
\newcommand{\SC}{\cC}
\newcommand{\ST}{\cT}
\newcommand{\rd}{{\rm d}}
\newcommand{\eins}{{\bf 1}}
\renewcommand{\I}{\eins}
\renewcommand{\SB}{\cB}
\newcommand{\gl}{\lambda}
\newcommand{\SR}{\cmc{R}}
\newcommand{\cst}{\star}
\newcommand{\ri}{{\rm i}}
\newcommand{\SU}{\cmc{U}}
\newcommand{\SI}{\cmc{I}}

\usepackage[varg]{txfonts}

\newcommand{\pn}{\par \noindent}

\begin{document}




\title{Transformations of Lévy Processes}

\author{Michael Schürmann, Michael Skeide\thanks{MS is supported by research funds of the Italian MIUR and of the University of Molise.}
\
and Silvia Volkwardt}

\date{November 2007, revised January 2008}

{
\renewcommand{\baselinestretch}{1}
\maketitle



\begin{abstract}\noindent
A Lévy process on a 
\nbd{*}bialgebra is given by its generator, a conditionally positive hermitian linear functional
vanishing at the unit element.
A 
\nbd{*}algebra homomorphism $\vk$ from
 a 
\nbd{*}bialgebra $\SC$ to a 
\nbd{*}bialgebra $\SB$ with the property that $\vk$ respects the counits maps generators
on $\SB$ to generators on $\SC$. A tranformation between the corrresponding two Lévy processes is given by
forming infinitesimal convolution products. 
This general result is applied to various situations, e.g., to a 
\nbd{*}bialgebra and its associated primitive tensor
\nbd{*}bialgebra (called \lq generator process\rq ) as well as its associated group-like 
\nbd{*}bialgebra (called Weyl-
\nbd{*}bialgebra).
It follows that a Lévy process on a 
\nbd{*}bialgebra can be realized on Boson Fock space as the infinitesimal convolution
product of its generator process such that the vacuum vector is cyclic for the Lévy process.
Moreover, we obtain convolution approximations of the Azéma martingale by the Wiener process and vice versa.

\end{abstract}


}


%

\section{Introduction}

A stochastic process $X_t : E \to G$, $t \geq 0$,
over some probability space $E$ taking values in a (topological)
group $G$ is called a (stationary) L\'evy process on $G$ if the increments
$X_{st} = X_s ^{-1} X_t$, $0 \leq s \leq t$, of disjoint intervals $[s,t)$
are independent, if the distribution of $X_{st}$ only depends on $t - s$ (stationarity), and if,
for $t \to 0+$, we have that $X_t$ converges in law to the Dirac measure $\delta _e$ concentrated 
at the unit element $e \in G$.
From an algebraic point of view this can immediately be generalized to stochastic
processes $(X_{st})_{0 \leq s \leq t}$ taking values in a monoid $G$ where the additional
evolution equation
$X_{rs} X_{st} = X_{rt}$ is postulated.
These \lq classical\rq \ L\'evy processes are commutative in the following sense.
If we replace $G$ and $E$ by  suitable 
\nbd{*}algebras of functions (on $G$ and $E$; e. g. replace
$G$ by $\RL ^{\infty} (G)$ and
$E$ by $\RL ^{\infty} (E )$)
then $X_{st} : E \to G$ will give a 
\nbd{*}algebra homomorphism
mapping a function $f$ on $G$ to the function $f \circ X_{st}$ on $E$.
The $j_{st}$ form a commutative process because they are defined on a commutative 
\nbd{*}algebra.
Replacing the monoid $G$ by a \nbd{*}bialgebra and the classical probability space $E$ by what is called
a quantum probability space, the notion of a quantum L\'evy process (QLP) on a \nbd{*}bialgebra over a
quantum probability space can be introduced; cf. \cite{ASchvW}.

The representation theorem for such processes \cite[Theorem 2.5.3]{MSchue93} says that they can always be
realized on a Boson Fock space as solutions to quantum stochastic differential equations 
in the sense of Hudson and Parthasarathy \cite{HuPa84}.
As pointed out in \cite{Ske05b} QLPs can also be viewed as tensor product systems of type I in the sense of W. Arveson \cite{Arv}. They are (up to 
stochastic equivalence) uniquely determinded 
by their generators which are the hermitian, normalized conditionally positive linear functionals on the underlying
\nbd{*}bialgebra. 
In this paper we are mainly interested in the following situation.
If there are given two bialgebras and an algebra homomorphism between them with the additional property that
the homomorphism respects the counits, then generators are transformed into generators. 
The question arises how the two 
QLPs  given by the two generators can be transformed into each other.
Using infinitesimal convolution products,
we establish a transformation on the level of the QLPs.

We describe very briefly what we do in a slightly simplified setting. (For instance, the example about Azéma martingales in Section \ref{AzemaSEC} fits into that simplified setting. For a precise description of the general situation see Sections \ref{prelSEC} and \ref{resSEC}.) In this simplified setting the situation is as follows: Suppose $(\cB,\Delta,\delta)$ is a \nbd{*}bialgebra. Then the comultiplication $\Delta$ induces a convolution $\star$ for algebra-valued linear mappings on $\cB$; see Section \ref{prelSEC}. Among all the properties a QLP $j=\bfam{j_{s,t}}_{0\le s\le t<\infty}$ satisfies, there is also the equality
\beqn{
j_{s,t}(b)
=
j_{t_0,t_1}\star\cdots\star j_{t_{n-1},t_n}(b)
}\eeqn
for all $s=t_0<t_1<\ldots<t_{n-1}<t_n=t$. Suppose on $\cB$ there is a second comultiplication $\Delta'$. We shall show that, in the canonical representation of $j$ on a pre-Hilbert space $D$ with cyclic vector $\Om$, the expressions
\beqn{
j_{t_0,t_1}\star'\cdots\star'j_{t_{n-1},t_n}(b)\Om
}\eeqn
(with the convolution with respect to $\Delta$ replaced by the convolution with respect to $\Delta'$) form a Cauchy net over the partions of the interval $[s,t]$. From this it easy to show that their limits, which we denote by $k_{s,t}(b)\Om$ determine on the their linear hull a unique QLP $k$ over $(\cB,\Delta',\delta)$, the \it{tranform} of $j$. Moreover, we shall show that under suitable cyclicity conditions this procedure can be reversed. See Theorem \ref{transthm} for a precise formulation in a more general context.

The transformation has various applications.
For example, there are two QLPs associated with a given QLP in a natural way.
One is the QLP's Weyl operator type process, the other is the generator process of the QLP which is 
composed of creation, annihilation and preservation processes on Boson Fock space.
The Weyl type process can be used to show in a nice way why the result of
Skeide \cite{Fra06} holds which says that the vacuum vector is always cyclic for the QLP.
The generator process allows for a construction of the QLP as a product system by infinitesimal
 convolution products as a kind of
multiplicative stochastic integral.
Both types of processes admit direct realizations on the Boson Fock space. Writing down the backwards transformations provides two different new proofs of the fact that every QLP may be realized as a (cyclic) process on a Boson Fock space. The relation with the generator process even reproves the fact that the original process fulfills a quantum stochastic differential equation.
Another application is the approximation of the Az\'ema martingales by infinitesimal
convolution products of the Wiener process, and \it{vice versa}.

In Section \ref{prelSEC} we repeat the necessary definitions that, in Section \ref{resSEC}, are used to formulate the transformation theorem. Section \ref{resSEC} also provides the constructions of several related \nbd{*}bialgebras, necessary for the applications. Section \ref{transSEC} presents the proof of the transformation theorem, Section \ref{applSEC} its applications.

\section{Preliminaries}\label{prelSEC}

An \hl{involutive} or \hl{\nbd{*}vector space} is a vector space $V$ with an involution, $\rm{i.e.}$, an anti-linear
 mapping $v\mapsto v^*$ on $V$ satisfying $(v^*)^*=v$. A \hl{\nbd{*}algebra} is an algebra $\cA$ which is also
 a \nbd{*}vector space such that $(ab)^*=b^*a^*$ for all $a,b\in\cA$.
If $\cA$ is a \nbd{*}algebra, then so is $\cA\otimes\cA$ with involution defined by $(a_1\otimes a_2)^*=a_1^*\otimes a_2^*$.

A complex vector space $\cC$ is a \hl{coalgebra} if there are linear maps $\Delta\colon\cC\rightarrow\cC\otimes\cC$ and
 $\delta\colon\cC\rightarrow\C$, called the \hl{coproduct} and \hl{counit} respectively, satisfying 
\baln{
(\Delta\otimes id)\circ\Delta&=(id\otimes\Delta)\circ\Delta\quad\text{(coassociativity)}   \\
(\delta\otimes id)\circ\Delta&=id=(id\otimes\delta)\circ\Delta\quad\text{(counit property)}.
}
\ealn
Following Sweedler we frequently
use the notation $c_{(1)}\otimes c_{(2)}$ for $\Delta(c)$ surpressing both summation and indices. 
Let $\Delta_0 := \gd$, $\Delta_1:=I_\cC$, and for $n\ge2$ define
$$
\Delta_n = (\GD _{n - 1} \ot \id ) \circ \GD .
$$
Sweedler's notation extends to writing $c_{(1)}\otimes c_{(2)}\otimes\cdots \otimes c_{(n)}$ for $\Delta_n(c),\thinspace n\geq 1$. 

Sometimes we shall need to equip also the \hl{conjugate} vector space $\ol{\cC}$ with a coalgebra structure. Note that the canonical bijection $i=i_1\colon c\mapsto c$ from $\cC$ to $\ol{\cC}$ is an anti-linear isomorphism. The same is true for the canonical bijections $i_n$ from the \nbd{n}fold tensor power of $\cC$ to the \nbd{n}fold tensor power of $\ol{\cC}$. Using $\ol{\ol{\cC}}=\cC$ we shall write $i_n^{-1}=i_n$. Note that $i_n\otimes i_m=i_{n+m}$ (where the tensor product of antilinear mappings is well-defined). By $i_0$ we denote complex conjugation of $\C$. It is, then, easy to convince oneself that $\overline{\delta}:=i_0\circ\delta\circ i_1$ and $\overline{\Delta}:=i_2\circ\delta\circ i_1$ make $(\ol{\cC},\overline{\Delta},\overline{\delta})$ a coalgebra. We shall use the notation $\ol{c}=i_1(c)$, so that $\ol{c_1\otimes\ldots\otimes c_n}=\ol{c}_1\otimes\ldots\otimes\ol{c}_n$

We shall also need the \hl{tensor product} $(\cC_1\otimes\cC_2,\Delta,\delta)$ of two coalgebras $(\cC_1,\Delta_1,\delta_1)$ and $(\cC_2,\Delta_2,\delta_2)$, where $\delta:=\delta_1\otimes\delta_2$ and $\Delta:=(\id\otimes\tau\otimes\id)\circ(\Delta_1\otimes\Delta_2)$ and $\tau$ denotes the flip $c\otimes d\mapsto d\otimes c$.

A \hl{\nbd{*}bialgebra} $(\cB,\Delta,\delta)$ is a coalgebra which is also a unital \nbd{*}algebra, and in such a way that
 $\Delta$ and $\delta$ are \nbd{*}algebra homomorphisms.
If $\cA$ is a unital \nbd{*}algebra with the multiplication map $M\colon\cA\otimes\cA\rightarrow\cA$ 
defined by setting $M(a_1\otimes a_2)=a_1a_2$, then we define the \hl{convolution} of two linear 
mappings $j,k\colon\cB\rightarrow\cA$ by $j\star k:=M\circ(j\otimes k)\circ\Delta$. 
In particular, 
the convolution of two linear functionals $\vp$ and $\psi$ on $\cB$ is $\vp\star\psi=(\vp\otimes\psi)\circ\Delta$.
 Unitality for a bialgebra $(\cB,\Delta,\delta)$ means that it is unital as an algebra, $\rm{i.e.}$, there exists $\bf{1}\in\cB$ 
such that $M(b\otimes\bf{1})=M(\bf{1}\otimes b)=b$ for all $b\in\cB$ and the coproduct and counit are unital, 
$\rm{i.e.}$, $\Delta(\bf{1})=\bf{1}\otimes\bf{1}$ and $\delta(\bf{1})=1$. 
We only consider unital algebras.

Let $(\cA,\Phi)$ be a \hl{quantum probability space}, 
that is, a unital \nbd{*}algebra with a \hl{state} 
 (a normalized positive linear functional $\Phi\colon\cA\rightarrow\C$).
 A \hl{quantum stochastic process} $j=\bfam{j_i}_{i\in I}$, indexed by some index set $I$, 
is a family of \hl{quantum random variables} $j_i$ (that is, of unital \nbd{*}algebra homomorphisms
 $j_i\colon\cB\rightarrow\cA$). By $\vp_i:=\Phi\circ j_i$ we denote the \hl{distribution} of $j_i$. 
The notion of independence used for quantum Lévy processes on \nbd{*}bialgebras in this paper is the tensor independence.
 A \hl{stationary quantum Lévy process} on $\cB$ over $\cA$ is a quantum stochastic process 
$j=\bfam{j_{s,t}}_{0\le s\le t<\infty}$, satisfying the following four conditions.
\begin{itemize}
\item[]
\begin{itemize}
\item[(LP1)]
The increments $j_{s,t}$ of disjoint intervals $\LO{s,t}$ are \hl{tensor independent} in $\Phi$, that is,
\baln{
\Phi\bfam{j_{s_1,t_1}(b_1)\cdots j_{s_n,t_n}(b_n)}
~&=~
\vp_{s_1,t_1}(b_1)\cdots\vp_{s_n,t_n}(b_n)\medspace\text{for all}\medspace n\in\N, b_k\in\cB\medspace\text{and}\\
[j_{s_k,t_k}(b_1),j_{s_l,t_l}(b_2)]&=0\medspace\text{for all}\medspace k\neq l\medspace\text{and all}\medspace b_1,b_2\in\cB,
}\ealn
whenever $k\ne\ell\Rightarrow\LO{s_k,t_k}\cap\LO{s_\ell,t_\ell}=\emptyset$.

\item[(LP2)]
The increments are \hl{stationary}, that is, $\vp_{s,t}=\vp_{0,t-s}$ for all $0\le s \le t $.

\item[(LP3)]
The process is \hl{continuous} in $\Phi$, that is, $\lim_{t\to0}\vp_{0,t}(b)=\delta (b)$ for all $b \in \SB$.

\item[(LP4)]
The $j_{s,t}$ are increments under convolution, 
that is, $j_{r,s}\star j_{s,t}=j_{r,t}$ for all $0\le r \le s \le t$ 
and $j_{t,t}(b)=\delta(b)\U$ for all $0\le t<\infty$.
\end{itemize}
\end{itemize}

\lf\noindent
In the sequel, for a stationary quantum Lévy process we will simply say \hl{Lévy process}.
We observe that by (LP1) and (LP4) every Lévy process fullfills the condition:
\begin{itemize}
\item[]
\begin{itemize}
\item[(LP4')]
$\vp_{r,s}\star\vp_{s,t}=\vp_{r,t}$ for all $0 \le r \le s \le t$ and $\vp_{t,t}=\delta$.
\end{itemize}
\end{itemize}
Therefore, by (LP2) and (LP3) the states $\vp_t:=\vp_{0,t}$ form 
a weakly continuous semigroup under convolution. By (LP1), (LP2) and (LP4)
 this convolution semigroup determines all \hl{joint moments} (that is exactly
 all expressions of the form of the left-hand side of the first equation of (LP1),
even if we drop the condition that the $( s_k , t_k ]$ are mutually disjoint). In other words,
 two Lévy processes are \hl{stochastically equivalent}, if and only if they have
 the same convolution semigroup. 
We can associate a \hl{generator} $\psi$ with a 
convolution semigroup through $\varphi_t=e_\star^{t\psi}$ for all $t\geq0$. Then $\psi$ is a 
linear functional on $\cB$, satisfying $\psi(\bf{1})=0$, and it is conditionally positive and hermitian. 
Thus, Lévy processes on \nbd{*}bialgebras can also be characterized (up to equivalence) by their generator. 

Let $D$ be a pre-Hilbert space and denote by $\sL^a(D)$ the \nbd{*}algebra of adjointable operators on $D$.
 If $\Om$ is a unit vector in $D$, then  $\bfam{\sL^a(D),\AB{\Om,\cdot \Om}}$ is a quantum probability space. 
We we call it a \hl{concrete} quantum probability space and write it as $(D,\Om)$. 
If a Lévy process $j$ takes values in a concrete quantum probability space, then we say $j$ is a \hl{concrete} Lévy process.
 By GNS-construction every quantum probability space $(\cA,\Phi)$ gives rise to a concrete quantum probability space $(D,\Om)$,
 determined uniquely by the properties that there is a \nbd{*}representation 
$\pi\colon\cA\rightarrow\sL^a(D)$ such that $\Phi=\AB{\Om,\pi(\cdot )\Om}$ and that $\Om$ is cyclic for $\cA$, that is,
 $\pi(\cA)\Om=D$. 
Consequently, every Lévy process gives rise to a concrete Lévy process \hl{over} $(D,\Om)$. 
In these notes we will consider only concrete Lévy processes and we will leave out the word concrete. 
We will say the Lévy process is \hl{cyclic}, if $\Om$ is cyclic for the \nbd{*}subalgebra
\beqn{
\cA_j
~:=~
\ls\BCB{j_{t_0,t_1}(b_1)\cdots j_{t_{n-1},t_n}(b_n)\colon n\in\N, b_k\in\cB, 0=t_0 \le \cdots \le t_n}
}\eeqn
of $\sL^a(D)$. (Recall that $j_{t,t}(b)=\delta(b)\U$. So, the case $t_{k-1}=t_k$
 can be excluded. Also, the case with $t_0>1$ can easily be achieved by putting $b_1=\U$.)
Notice that by (LP1) this space does not change, if we allow that the disjoint intervals are not consecutive.
 By restricting to the invariant subspace $\cA_j\Om$ of $D$ that is generated by the process from $\Om$, 
we obtain from every Lévy process over $D$ a cyclic Lévy process on $\cA_j\Om = D_j$.

By a GNS-type construction applied 
to a generator $\psi$ on $\cB$ we obtain a pre-Hilbert space $K$,
 a surjective mapping $\eta\colon\cB\rightarrow K$ and a \nbd{*}representation
 $\rho\colon\cB\rightarrow\sL^a(K)$ such that 
\beqn{
\eta(ab)=\rho(a)\eta(b)+\eta(a)\delta(b)
}\eeqn
and
\beq{
\label{korand}
-\langle\eta(a^*),\eta(b)\rangle=\delta(a)\psi(b)-\psi(ab)+\psi(a)\delta(b)
}\eeq
for all $a,b\in\cB$. The specified triple $(\rho,\eta,\psi)$ is called a \hl{surjective Lévy triple}. 
There is a one-to-one correspondence between 
Lévy processes (modulo equivalence) on $\cB$, convolution semigroups of states on $\cB$, generators on $\cB$
 and surjective Lévy triples on $\cB$ 
(modulo unitary equivalence).

Of course, for every convolution semigroup $\vp=\bfam{\vp_t}_{t\in\R_+}$ there is (up to unitary equivalence) 
at most one cyclic Lévy process. 
(Unitary equivalence is much stronger than stochastic equivalence.) Effectively, if $j$ is a cyclic process on 
$(D,\Om)$ which fulfills (LP1) - (LP3) and (LP4'), then is not difficult to show that also (LP4) holds. 
By a GNS-type construction Schürmann \cite[Proposition 1.9.5]{MSchue93} shows that every 
convolution semigroup of states on a \nbd{*}bialgebra there is a (unique up to unitary equivalence) cyclic Lévy process (even whithout continuity). 
This construction involves the GNS-construction of all $\vp_t$, their tensor products and an inductive limit over the interval 
partitions of $\R_+$. However, it is completely algebraic and does not involve analytic tools.
On the contrary, \cite[Theorem 2.5.3]{MSchue93} constructs 
a Lévy process on a (symmetric) Fock space $\G(L^2(\R_+,K))$ as the
 solution of a (quite an involved
 system of) quantum stochastic differential equation(s) in the sense of 
Hudson and Parthasarathy \cite{HuPa84}. For quite a long time it was an open problem,
 to decide whether Fock space and differential equation can be set in such a way that the 
Fock vacuum is cyclic for the resulting Lévy process. Only quite recently and simultaneously, Franz,
 Schürmann and Skeide came up, not with just one, but with a whole bunch of proofs for the affirmative answer.

The proof due to Skeide (see Franz \cite[Theorem 1.21]{Fra06}) uses in an essential way the representation on 
the Fock space and the differential equation of
 \cite[Theorem 2.5.3]{MSchue93} and shows that for every $b\in\cB$ with $\delta(b)=1$ the vectors
\beq{\label{expappr}
j_{t_0,t_1}(b)\cdots j_{t_{n-1},t_n}(b)\Om  ,
}\eeq
$s=t_0 \le t_1 \le \ldots \le t_{n-1} \le t_n=t$, converge over the interval 
partitions of $\LO{s,t}$ to an exponential vector of the form $\exp(k\I_{\LO{s,t}})$ 
where $k\in K$ is a vector depending on $b$. 
(Cyclicity is, then, a simply consequence of Skeide's proof in \cite{Ske00a} of a result due to Parthasarathy and Sunder \cite{PaSu98}.)
Immediately, from this construction,  the idea emerged to construct an explicit isomorphism from 
the space of the abstract Lévy process of \cite[Proposition 1.9.5]{MSchue93} to the Fock space of the Lévy process 
obtained via \cite[Theorem 2.5.3]{MSchue93}. Namely, if in \eqref{expappr} we replace $j$ and $\Om$ with the
 abstract process $j'$ and its cyclic vector $\Om'$, we know from \cite[Theorem 1.21]{Fra06} that they converge. 
Sending the limit to $\exp(k\I_{\LO{s,t}})$ establishes a unitary from the abstract representation space $D'$ to the Fock space. 
If we can manage to do this without using \cite[Theorem 1.21]{Fra06}, then we will obtain a direct proof of representability of
 the Lévy process as cyclic process on the Fock space. 

The idea for a transformation of a (cyclic) Lévy process originates in the following observation.
Let us denote by $\cB_1:=\CB{b\in\cB\colon\delta(b)=1}$ the set of all elements in $\cB$ to which \eqref{expappr} applies. 
Suppose the element $b\in\cB_1$ is \hl{group-like}, that is, $\Delta(b)=b\otimes b$. 
(Note that $b\in\cB$ being group-like, the counit property forces $b=0$ or $b\in\cB_1$.) Then
\beqn{
j_{t_0,t_1}(b)\cdots j_{t_{n-1},t_n}(b) =
j_{t_0,t_1}\star\cdots\star j_{t_{n-1},t_n}(b) =
j_{s,t}(b)
}\eeqn
so that the limit is over a constant and gives back what $j_{s,t}(b)$ does to the cyclic vector.
 In general, there need not be group like elements in $\cB_1$, and if, then they need not generate $\cB$. However, 
if we were able to define a different comultiplication on $\cB$ for which all elements in $\cB_1$ are group-like, then
\beqn{
k_{s,t}(b)\Om =
\lim j_{t_0,t_1}(b)\cdots j_{t_{n-1},t_n}(b)\Om
}\eeqn
would define a family of homomorphisms $k_{s,t}$ that form a Lévy process with respect to the group-like comultiplication. 
In other words, we transformed one Lévy process into another.

It is easy to give a direct realization of such a group-like process on a suitable Fock space; see Section 4.1. 
Thus, provided 
that the process $k$ acts cyclic on $\Om$, we would find the representation theorem. 
The easiest way to establish cyclicity is 
to reconstruct $j$ from $k$ by a reverse transformation.
 Recall that the construction of $k$ involved replacing the original comultiplication 
with one that makes all $b\in\cB_1$ into group-like elements so that $j_{t_0,t_1}(b)\cdots j_{t_{n-1},t_n}(b)$ is 
nothing but
 $j_{t_0,t_1}\star'\cdots\star'j_{t_{n-1},t_n}$ with respect to the new comultiplication. 
Now we do just the opposite and look at the limit of
\beq{\label{invappr}
k_{t_0,t_1}\star\cdots\star k_{t_{n-1},t_n}(b)\Om
}\eeq
for the original comultiplication. If this reverse transformation gives back $j$,
then, knowing that the representation space of the 
intermediate group-like process $k$ is isomorphic to a Fock space, we will know that
 also the representation space of $j$ is a Fock space.
Technically, in general, it is not possible to equip $\cB$ directly with a comultiplication
that makes the elements of $\cB_1$ group-like.
 However, it is possible to associate with every \nbd{*}bialgebra 
$\cB$ its \hl{group-like \nbd{*}bialgebra} $\C\cB_1$. 
The vector space $\C\cB_1$
 contains the set $\cB_1$ as a basis consisting entirely of group like elements.
 And the $k_{s,t}(b)\Om$ defined on elements of $\cB_1$ determine
 a unique Lévy process on $\C\cB_1$. But now the $k_{st}$ do no longer define a linear 
mapping $\cB\rightarrow\sL^a(D)$. 
(They do define a linear mapping $\C\cB_1\rightarrow\sL^a(D')$ where $D'$ is the linear span in $\ol{D}$ of
 what the $k_{s,t}(b)$ generate from 
$\Om$.) 
So the convolutions in \eqref{invappr} with respect to the comultiplication of $\cB$ do no longer have a meaning.
 The problem is solved if
 we associate again with $\cB$ a special kind of \nbd{*}bialgebra; see example \ref{extarev}. We will equip this 
tensor \nbd{*}bialgebra with a certain
 comultiplication, so that the convolutions in \eqref{invappr} are defined with respect to this comultiplication.

\section{Statement of results}\label{resSEC}

\lf\noindent
We start our considerations with a cyclic Lévy process on $(\cB,\Delta,\delta)$ whose generator is $\psi$. 
Furthermore, there are given another \nbd{*}bialgebra $(\cC,\Lambda,\lambda)$ and 
a unital \nbd{*}algebra homomorphism $\vk\colon\cC\rightarrow\cB$ which \hl{preserves the counit},
 $\rm{i.e.}$, $\delta\circ\vk=\lambda$. 
Since $\vk ({\bf 1}) = {\bf 1}$ it is easy to see that this last property is equivalent to the condition
$\vk (\cC _0 ) \subset \cB _0$ where $\SC _0 = \ker \gl$, $\SB _0 = \ker \gd$.
A generator of a Lévy process on $\cB$ can be lifted via $\vk$ to a generator $\psi\circ\vk$ of a Lévy process on $\cC$. 
Therefore, the question arises, 
what is the relationship between the two Lévy processes? We will show how the second process can be computed from the first
 one and vice versa.

\bex
{
\label{exta}
\hl{(Primitive tensor \nbd{*}bialgebra associated with a \nbd{*}bialgebra)}
\par\noindent
For a vector space $V$ the \hl{tensor algebra} $\sT(V)$ is the vector space 
\[\sT(V) =\bigoplus_{n\in\N} V^{\otimes n}\] 
where $V^{\otimes n}$ denotes the $n$-fold tensor product of V with itself, $V^{\otimes 0}=\C$, 
with the multiplication given by $(v_1\otimes\cdots \otimes v_n,v_{n+1}\otimes\cdots \otimes 
v_{n+r})\mapsto v_1\otimes\cdots \otimes v_n\otimes v_{n+1}\otimes\cdots \otimes v_{n+r}$ for $n,r\in\N, v_1, \dots, 
v_n, v_{n+1}, \dots, v_{n+r}\in V$. The tensor algebra satisfies the following universal property. 
There exists an embedding $\iota\colon V\rightarrow\sT(V)$ of $V$ to $\sT(V)$ such that any linear mapping $f$ from $V$ 
into an algebra $\cA$ can be uniquely extended to an \emph{algebra homomorphism} $\sT(f)\colon\sT(V)\rightarrow\cA$ such 
that $\sT(f)\circ\iota(v)=f(v)$ for all $v\in V$. 
Conversely, any algebra homomorphismus $g\colon\sT(V)\rightarrow\cA$ is 
uniquely determined by its restriction to $V$. 
In a similar way, an involution on $V$ gives rise to a unique extension as an involution on $\sT(V)$. 
Thus, for a \nbd{*}vector space V we can form the tensor \nbd{*}algebra $\sT(V)$.
This can be used to define a unique \nbd{*}bialgebra structure on $\sT(V)$ such that all elements in $V$ are primitive, 
$\rm{i.e.}$, the extended mappings $\Lambda\colon V\rightarrow\sT(V)\otimes\sT(V), v\mapsto v\otimes\bf{1}+\bf{1}\otimes v$ and 
$\lambda\colon V\rightarrow \cC, v\mapsto 0$ define the comultiplication and the counit on $\sT(V)$.

Let $(\cB,\Delta,\delta)$ be any \nbd{*}bialgebra. The set $\cB_0=\{b\in\cB\colon\delta(b)=0\}$ is an \nbd{*}ideal of $\cB$. 
The tensor 
\nbd{*}bialgebra $\sT(\cB_0)$ 
over $\cB _0 )$ is a \nbd{*}bialgebra with the above comultiplication and counit. 
So the second \nbd{*}bialgebra $\cC$ is  $(\sT(\cB_0),\Lambda,\lambda)$ which is called the \hl{primitive 
tensor \nbd{*}bialgebra} associated with $\cB$. The counit preserving \nbd{*}algebra homomorphism $\vk$ is defined by 
$\vk(b_1\otimes \cdots \otimes b_n)=b_1\cdots b_n$ for $b_1,\ldots , b_n\in\cB_0$.
}
\eex

\bex
\par\noindent
{
\label{extarev}
\hl{(Induced tensor \nbd{*}bialgebra associated with a \nbd{*}bialgebra)}
\par\noindent
Let $(\cB,\Delta,\delta)$ and $(\sT(\cB_0),\Lambda,\lambda)$ be the \nbd{*}bialgebras as in example \ref{exta}.
We can define another 
coalgebra structure on $\sT(\cB_0)$.
Denote by
\beqn{
E : \SB _0 \oplus \SB _0 \oplus (\SB _0 \ot \SB _0 ) \to  \sT (\SB _0 ) \ot \sT (\SB_0 )
}\eeqn
the canonical embedding coming from the identification of
$\SB _0$ with $\SB _0 \ot \eins$ and 
$\eins \ot \SB _0$ respectively and $\SB _0 \ot \SB _0 \subset
\sT (\SB _0 ) \ot \sT (\SB _0 )$.
Moreover, consider the restriction $\GD _0$ of $\GD$ to $\SB _0$. 
Then 
\beqn{
\GD _0 : \SB _0 \to \SB _0 \oplus \SB _0 \oplus (\SB _0  \ot \SB _0  )
}\eeqn
and
$(\sT (\SB _0 ) , \sT ( E \circ \GD _0 ) , \sT (0))$ is a \nbd{*}bialgebra.
We can understand this \nbd{*}bialgebra as  a \lq big version\rq \ of $\cB$
 and so $(\sT(\cB_0),\sT(\Delta _0 ),\sT(0))$ is called the \hl{induced tensor \nbd{*}bialgebra} associated with $\cB$.
In the context of the algebraic set-up the first \nbd{*}bialgebra is $(\sT(\cB_0),\Lambda, \sT (0))$ 
and the second \nbd{*}bialgebra is $(\sT(\cB_0),\sT(\Delta _0 ),\sT(0))$. The identity on $\sT(\cB_0)$ is an example
 of a counit preserving \nbd{*}algebra homomorphism $\vk$.
}
\eex

\bex{
\label{reversion}
\hl{(Reversion of the transformation)}
\pn
The reverse transformation of a Lévy process on $(\cC,\Lambda,\lambda)$ into a Lévy process on
 $(\cB,\Delta,\delta)$ requires a counit preserving \nbd{*}algebra homomorphism $\widetilde{\vk}$ 
which, roughly speaking, is the inverse of $\vk$. The construction of $\widetilde{\vk}$ assumes in addition the
 surjectivity of $\vk$. This implies 
$\vk (\cC _0 ) = \cB _0$ and
the existence of an injective linear 
\nbd{*}mapping
\[\upsilon\colon\cB _0 \rightarrow\cC_0 \quad\text{such that}\quad\vk\circ\upsilon=id_\cB.\]
The linear 
\nbd{*}mapping $\upsilon$ is not unique. Its existence follows from the existence of a self-adjoint
basis
$\bfam{b_i}_{i\in I}$ of the 
\nbd{*}vector space $\cB _0$, $I$ some index set.
 Choose $c_i\in\cC$ self-adjoint such that $\vk(c_i)=b_i$. This is possible since $\vk$ is surjectiv. 
Define the linear \nbd{*}map $\upsilon$ by $\upsilon(b_i)=c_i$. 
In view of the universal property of tensor algebras we extend the linear 
\nbd{*}map $\upsilon$ to a \nbd{*}algebra homomorphism 
\[\wt{\vk}=\sT(\upsilon )\colon\sT(\cB_0)\rightarrow\cC\]
to the induced tensor \nbd{*}bialgebra $\sT(\cB_0)$. The coalgebra structure on $\sT(\cB_0)$ is
 defined as in Example \ref{extarev} by
$\sT(\Delta)(b)=\Delta(b)$ and $\sT(\delta)(b)=\delta(b)$ for $b\in\cB$. Indeed, the \nbd{*}algebra 
homomorphism $\wt{\vk}$ preserves the counits. It is sufficient to show this for the generators of
 $\sT(\cB_0)$. For all $b\in\cB_0$ we have
\[\lambda\circ\upsilon(b)=\delta\circ\vk\circ\upsilon(b)=\delta\circ id_{\cB_0}(b)=0=\delta(b).\]
The above situation is described by
\[(\sT(\cB_0),\sT(\Delta),\sT(\delta))\medspace\xrightarrow[]{\wt{\vk}}\medspace(\cC,\Lambda,\lambda)
\medspace\xrightarrow[]{\vk}\medspace(\cB,\Delta,\delta).\]
}\eex

\bex
{
\label{exgl}
\hl{(Group-like \nbd{*}bialgebras)}
\par\noindent
For a set $S$ the \hl{vector space generated by $S$} is the vector space
\[\C(S):=\BCB{f\colon S\rightarrow\C : f(m)=0\medspace\text{for all but finitely many}\medspace s\in S}.\]
Assume in addition that $S$ is a monoid with identity $e\in S$. 
Since $S$ is a basis, the multiplication map $S\times S\rightarrow S$ induces a map $M\colon\C(S)\otimes\C(S)\rightarrow\C(S)$ that turns $\C(S)$ into an algebra with identity element $e\in S\subset\C(S)$.
Since $S$ is a basis of $\C(S)$ the mapping $M$ induces an algebra structure on $\C(S)$ 
with unit element $e$. The vector space generated by a set satisfies the following universal property. 
There exists an embedding $\iota\colon S\rightarrow\C(S)$ such that any mapping $\phi$ from $S$ to some vector space $Z$
 can be uniquely extended to a linear mapping $\ol{\phi}\colon\C(S)\rightarrow Z$ such that $\phi=\ol{\phi}\circ\iota$. 
This can be used to define a coalgebra structure on $\C(S)$.
We understand $S$ as a set of group-like elements.
 We extend the mappings $\Delta\colon S\rightarrow\C(S)\otimes\C(S),\thinspace\Delta(s)=s\otimes s$ and
$\delta\colon S\rightarrow\C,\thinspace\delta(s)=1$ to linear mappings on $\C(S)$. We will denote the comultiplication and the
 counit on $\C(S)$ again by $\Delta$ and $\delta$. Indeed, $\Delta$ and $\delta$ are algebra homomorphism since 
$\Delta(xy)=xy\otimes xy=(x\otimes x)(y \otimes y)=\Delta(x)\Delta(y)$ and $\delta(xy)=1=\delta(x)\delta(y)$ for all $x,y\in S$. 
An involution on $S$ can also be uniquely extended to an involution on $\C(S)$. Thus, for a \nbd{*}monoid $S$ we can
 form the \hl{group-like \nbd{*}bialgebra} $(\C(S),\Delta,\delta)$ over $S$.

Let $(\cB,\Delta,\delta)$ be a \nbd{*}bialgebra. The set $\cB_1=\{b\in\cB\colon\delta(b)=1\}$ is a 
\nbd{*}monoid with multiplication and involution of the 
\nbd{*}algebra $\cB$. 
Hence, $(\C(\cB_1),\Lambda,\lambda)$ is a \nbd{*}bialgebra, the so called \hl{group-like \nbd{*}bialgebra
 associated to} $(\cB,\Delta,\delta)$. 
In the sequel, we write $\widehat{b}$ for the element $b$ in $\cB_1\subset\C\cB_1$. The comultiplication $\Lambda$ 
and the counit $\lambda$ on 
$\C(\cB_1)$ are defined by $\Lambda(\thinspace\wh{b}\thinspace)=\wh{b}\otimes\wh{b}$ and
 $\lambda(\thinspace\wh{b}\thinspace)=1$ for $\wh{b}\in\C(\cB_1)$. $\cB_1$ is equal to the set of 
all group-like elements in
 $\C\cB_1$, $\rm{i.e.}$, $\cB_1=\{0\neq\wh{b}\in\C(B_1)\colon\Lambda(\thinspace\wh{b}\thinspace)=\wh{b}\otimes\wh{b}\}.$
 Therefore, we have
\[(\sT(\cB_0),\sT(\Delta),\sT(\delta))\medspace\xrightarrow[]{\wt{\vk}}\medspace(\C(\cB_1),\Lambda,\lambda)\medspace\xrightarrow[]
{\vk}\medspace(\cB,\Delta,\delta)\]
where the counit preserving \nbd{*}algebra homomorphism $\vk$ and $\wt{\vk}$ are defined by $\vk(\thinspace\widehat{b}\thinspace)=b$ for 
$b\in\cB_1$ and $\wt{\vk}(b)=\wh{b+\bf{1}}-\wh{\bf{1}}$ for $b\in\cB_0$. 
Now we are able to express the reverse transformation \ref{invappr} by 
$(k_{t_0,t_1} \circ \widetilde{\vk}) \star_{_{\sT (\Delta )}} \cdots \star_{_{\sT (\Delta )}} 
(k_{t_{n-1},t_{n}} \circ \widetilde{\vk}) (b) \Omega$ for $b\in\cB_0$.
}
\eex

\lf
In the sequel, $\eZ_{st}$ denotes the set of all \hl{partitions} of an interval $[s,t]\subset\R^+$. 
Let $\alpha=\{s=t_0<t_1<\dotsb<t_{n-1}<t_n=t\}$ be a partition of $[s,t]$ and define 
$$
\norm{\alpha}=\max\{t_{j+1}-t_j\mid 0\leq j \leq n-1\}.
$$ 
We turn $\eZ_{st}$ into a directed set by writing 
$
\alpha_1\prec \alpha_2:\Leftrightarrow\alpha_1\subset\alpha_2.
$

\bthm\label{transthm}
Let $(\cB,\Delta,\delta)$ be a \nbd{*}bialgebra and let $(j_{s,t})_{0\leq s\leq t}$ be the unique cyclic Lévy
 process over $(D_j ,\Om)$ whose convolution semigroup is given by a generator $\psi$. Let $(\cC,\Lambda,\lambda)$ 
be another \nbd{*}bialgebra and let $\vk\colon\cC\rightarrow\cB$ be a unital \nbd{*}algebra homomorphism which 
\hl{preserves the counit}, that is, $\delta\circ\vk=\lambda$. Denote by $H_k$ the Hilbert subspace of $\ol{D_j}$ defined by
\bmun{
H_k
~:=~
\cls\BCB{(j_{t_0,t_1}\circ\vk)(c_1)\cdots (j_{t_{n-1},t_n}\circ\vk )(c_n )\Om\colon
\\
n\in\N,c_1,\dotsc,c_n\in\cC,0\le s \le t, s=t_0\leq t_1\leq\dotsb \leq t_{n-1}\leq t_n=t}.
}\emun
\begin{enumerate}
\item\label{trans1}
For every $c\in\cC$ and $0\le s\le t$ the net $\bfam{\vt_\alpha(c)}_{\alpha\in\eZ_{st}}$
converges in norm to an element in $H_k$ where
\beq{\label{vtdef}
\vt_\alpha(c)=(j_{t_0,t_1}\circ\vk)\star\cdots \star(j_{t_{n-1},t_n}\circ\vk)(c)\Om.
}\eeq
 
Moreover, setting
\beqn{
k_{s,t}(c)\Om
~:=~
\lim_\alpha 
\vartheta _{\alpha} \, (c)
}\eeqn
determines a unique cyclic Lévy process $k=\bfam{k_{s,t}}_{0\le s\le t<\infty}$ on $\cC$ over 
a dense subspace $(D_k,\Om)$ of $H_k$. The convolution semigroup of this process has generator $\psi \circ\vk$.

\item\label{trans2}
Let $\bfam{k_{s,t}}_{0\leq s\leq t }$ be the cyclic Lévy process on $(\cC,\Lambda,\lambda)$ 
over $(D_k,\Om)$ as constructed in the first part of the theorem. 
Assume in addition that $\vk$ is surjective.
 Let $(\sT(\cB_0),\sT(\Delta),\sT(\delta))$ be the induced tensor 
\nbd{*}bialgebra associated with $(\cB,\Delta,\delta)$ and 
let $\widetilde{\vk}\colon\sT(\cB_0)\rightarrow\cC$ be like in Example \ref{reversion}.

For every $b\in\cB$ and $0\le s \le t$ the net $\bfam{\zeta_\alpha}_{\alpha\in\eZ_{st}}$ converges in norm to $j_{s,t}(b) \Om$ where
\[\zeta_\alpha (b) :=(k_{t_0,t_1} \circ \widetilde{\vk}) \star_{_{\sT (\Delta )}} 
\cdots \star_{_{\sT (\Delta )}} (k_{t_{n-1},t_{n}} \circ \widetilde{\vk}) (b) \Omega\] 
and $\bfam{j_{s,t}}_{0\leq s\leq t<\infty}$ is the original Lévy process on $(\cB,\Delta,\delta)$.
Moreover, we have $H_k = \ol{D_j }$.
\end{enumerate}
\ethm

\section{Proof of Theorem \ref{transthm}}\label{transSEC}

In principle, Part \ref{trans1} of Theorem \ref{transthm} is proved (and Part \ref{trans2} almost) if we show that the nets in \eqref{vtdef} are Cauchy. To that goal in Section \ref{lemSS} we prove a lemma about infinitesimal products in Banach algebras (an extension of ideas in \cite{LiSk05}) and a coalgebra version (appealing to the Fundamental Theorem of Coalgebras). These lemmas plus the algebraic Proposition \ref{BcCprop} allow to prove Proposition \ref{prop1}, which is the analytic heart of the proof of Theorem \ref{transthm}.

\subsection{Preparatory lemmas}\label{lemSS}

We start with a lemma that imitates, like in \cite{LiSk05}, proofs of the Trotter product formula.

\blem
\label{matrix1}
Let $\cA$ be a Banach algebra. Suppose we have a constant $R>0$ and a family $\bfam{A^{(\mu)}}_{\mu\in M}$ of functions
\[
r
~\longmapsto~
A_r^{(\mu)}=I+ r G + \eS _ r ^{(\mu)}
~\in~
\cA
\] 
on $\R_+$ where $G\in\cA$ and $\eS_r^{(\mu)}$ satisfies $\norm{\eS_r^{(\mu)}}\leq r^2 \frac{C^2}{2}$ for some constant $C$ not depending on $\mu\in M$ and all $r\le R$. Then for all intervals $[s,t]\subset\R_+$, all partitions $\alpha=\{s=t_0<t_1<\dotsb<t_{n-1}<t_n=t\}$ $(n\in\N)$ of $[s,t]$ with $\norm{\alpha}\le R$, and an arbitrary choice of elements $\mu_1,\ldots,\mu_n$ of $M$, we have
\beqn{
\norm{A_{t_1-t_0}^{(\mu_1)} \cdots A_{t_n-t_{n-1}}^{(\mu_n)}-e^{(t-s)G}}
~\le~
\norm{\alpha}(t-s)e^{(t-s) \max (\parallel G \parallel, C)}\frac{C^2+\norm{G}^2e^{\norm{\alpha}\,\norm{G}}}{2}.
}\eeqn
\elem

\proof
By assumption $\norm{A_r ^{(\mu_k)}}\leq 1 + r \norm{G} + r^2\frac{C^2}{2}\leq e^{r \max(\norm{G},C)}$, and thus
\[
\norm{A_{t_\ell-t_{\ell-1}}^{(\mu_\ell)} \cdots A_{t_k-t_{k-1}}^{(\mu_k)}}
\leq
e^{(t_k-t_{\ell-1})\max(\norm{G},C)}
\]
for all intervals $[s,t]\subset\R_+$, all partitions $\alpha_n$ of $[s,t]$, and all $1\le\ell<k\le n$. The next calculation (cf. \cite{LiSk05} proof of Proposition 3.3) is essential for the proof. We compute
\baln{
A_{t_1-t_0}^{(\mu_1)} \cdots A_{t_n-t_{n-1}}^{(\mu_n)}-e^{(t-s)G}
&
=
A_{t_1-t_0}^{(\mu_1)} \cdots A_{t_n-t_{n-1}}^{(\mu_n)}-e^{(t_1-t_0)G}\cdots e^{(t_n-t_{n-1})G}
\\
&
=
\sum_{j =1}^nA_{t_1-t_0}^{(\mu_1)} \cdots A_{t_{j-1}-t_{j-2}}^{(\mu_{j-1})}\Bfam{A_{t_j-t_{j-1}}^{(\mu_j)}-e^{(t_j-t_{j-1})G}}e^{(t_{j+1}-t_j)G}\cdots e^{(t_n-t_{n-1})G}.
}\ealn
We have
\baln{
\bnorm{A_{t_j-t_{j-1}}^{(\mu_j)}-e^{(t_j-t_{j-1})G}}
&
\leq 
\bnorm{A_{t_j-t_{j-1}}^{(\mu_j)} - I - (t_j-t_{j-1})G} + 
\bnorm{ I + (t_j-t_{j-1})G - e^{(t_j-t_{j-1})G}  }
\\
& \leq (t_j-t_{j-1})^2 \frac{C^2+\norm{G}^2e^{(t_j-t_{j-1})\norm{G}}}{2}. 
}\ealn
From this estimate, from the estimate preceding it, and from the estimate
\beqn{
\sum_{j=1}^n(t_j-t_{j-1})^2
\le
\norm{\alpha}\sum_{j=1}^n(t_j-t_{j-1})
=
\norm{\alpha}(t-s)
}\eeqn
the statement follows.\qed

\lf
There is a coalgebra version of Lemma \ref{matrix1} deduced from 
the \emph{Fundamental Theorem of Coalgebras}, which states that the coalgebra generated by a 
finite subset of a coalgebra is finite dimensional.
In the sequel, $\bf{L}(V,W)$ denotes the vector space of linear maps between vector spaces $V$ and $W$.
We put $\bf{L} (V, V) = \bf{L} (V)$.
Let $(\cC,\Delta,\delta)$ be a coalgebra and let $\psi\in\bf{L}(\cC,\C)$ be a linear functional on $\cC$. 
The map $T: \psi\mapsto(id\otimes\psi)\circ\Delta$ defines an
 injective unital algebra homomorphism from $(\bf{L}(\cC,\C),\star)$ to $(\bf{L}(\cC),\circ)$ with left inverse $\delta\circ \eins$.
 Moreover, each $T(\psi)$ leaves every sub-coalgebra of $\cC$ invariant. On an arbitrary finite-dimensional subcoalgebra $\cC_c\ni c$ of $\cC$ the the series $e^{M(\psi)}\restriction\cC_c:=\sum_{n=0}^\infty\frac{M(\psi)\restriction\cC_c}{n!}$ converges in any norm. By the Fundamental Theorem of Coalgebras for every $c\in\cC$ such a $\cC_c$ exists. We  deduce that the series 
\beq{
\label{convexp}
e_\star^\psi(c)
~:=~
\sum_{n=0}^\infty\frac{\psi^{\star n}}{n!}(c)
~=~
\delta\circ e^{T(\psi)}(c)
}\eeq
converges for all $\psi\in\bf{L}(\cC,\C)$ and all $c\in\cC$. Clearly, this limit of complex numbers cannot depend on the choice of $\cC_c$: see \cite{ASchvW}.

We now prove the coalgebra version of Lemma \ref{matrix1}.

\blem
\label{lem:lim}
Let $\cC$ be a coalgebra. Suppose we have a constant $R>0$ and a family $\bfam{f^{(\mu)}}_{\mu\in M}$ of functions
\[
r
~\longmapsto~
f_r^{(\mu)}=\delta+ r \psi + \eR _ r ^{(\mu)}
~\in~
\bf{L}(\cC,\C)
\] 
on $\R_+$ where $\psi\in\bf{L}(\cC,\C)$ and $\eR_r^{(\mu)}(c)$ satisfies $\abs{\eR_r^{(\mu)}(c)}\leq r^2 D_c$ for some constant $D_c>0$, depending on $c\in\cC$ but not on $\mu$, and all $r\le R$. Then there exist constants $C_c>0$ and $\Psi_c>0$ such that for all intervals $[s,t]\subset\R_+$, all partitions $\alpha_n=\{s=t_0<t_1<\dotsb<t_{n-1}<t_n=t\}$ $(n\in\N)$ of $[s,t]$ with $\norm{\alpha}\le R$, and an arbitrary choice of elements $\mu_1,\ldots,\mu_n$ of $M$, we have
\beqn{
\abs{f_{t_1-t_0}^{(\mu_1)}\star \cdots \star f_{t_n-t_{n-1}}^{(\mu_n)}(c)-e_\star^{(t-s)\psi}(c)}
~\le~
\norm{\alpha}(t-s)e^{(t-s) \max (\Psi_c, C_c)}\frac{C_c^2+\Psi_c^2e^{\norm{\alpha}\,\Psi_c}}{2}.
}\eeqn
\elem

\proof
Choose $b\in\cC$ and fix a finite-dimensional sub-coalgebra $\cC_b$ of $\cC$ containing $b$. 
Fix a norm on $\cC_b$. From the weak estimates $\abs{\eR_r^{(\mu)}(c)}\leq r^2 D_c$ we easily conclude the strong estimate $\norm{\eR_r^{(\mu)}}\leq r^2 D$ for a suitable constant $D$ for the linear functionals $\eR_r^{(\mu)}$ on $\cC_b$. (Just take your favorite elementary proof of the Uniform Boundedness Principle for finite-dimensional Banach spaces.) Consider the linear operator
\[A_r ^{(\mu)}:=T(f_r ^{(\mu)})\restriction\cC_b\]
on $\cC_b$, so
\[A_r ^{(\mu)}=I+r  G+\eS_r ^{(\mu)}\]
where $G:=T(\psi)\restriction\cC_b$ and $\eS_r ^{(\mu)}=T(\eR_r ^{(\mu)})\restriction\cC_b$.

$\bf{L}(\cC_b)$ is a Banach algebra with respect to the operator norm. Since $T$ is a bijection from $\bf{L}(\cC_b,\C)$ onto $T(\bf{L}(\cC_b,\C))\subset\bf{L}(\cC_b)$, and since all norms on finite-dimensional spaces are equivalent, $\eS_r ^{(\mu)}$ satisfies $\norm{\eS_r ^{(\mu)}}\leq r ^2 \frac{C^2}{2}$ for some constant $C$. In view of lemma \ref{matrix1} we obtain the claimed statement if we choose $C_c=C\sqrt{\norm{\delta}\norm{c}}$ and $\Psi_c=\norm{G}\sqrt{\norm{\delta}\norm{c}}$.\qed

\subsection{Proof of Part \ref{trans1} of Theorem \ref{transthm}}\label{p1SS}

Consider the Hilbert subspaces $(0 \leq s \leq t)$
\beqn{
H_{st} = \cls\BCB{j_{t_0,t_1}(b_1) \cdots j_{t_{n-1},t_n}(b_n) \Omega \mid n \in \mathbb{N}, \thinspace s = t_0 \leq t_1 \leq \cdots \leq t_n = t, 
\thinspace b_1, \dotsc ,b_n \in\cB}
}\eeqn
of $\ol{D_j}$ where $H_0 = \C$. 
Put $H_t = H_{0t}$.
Using the shift and the unit vector $\Om$, we define mappings $U_{st}\colon H_s \otimes H_t\rightarrow H_{s+t}$ by
\bmun{
U_{st}(j_{s_0,s_1}(b_1) \cdots j_{s_{n-1},s_n}(b_n) \Omega \thinspace \otimes \thinspace j_{t_0,t_1}(c_1) \cdots j_{t_{m-1},t_m}(c_m) \Omega ) \\
= j_{s_0,s_1}(b_1) \cdots j_{s_{n-1},s_n}(b_n) j_{t_0+s,t_1+s}(c_1) \cdots j_{t_{m-1}+s,t_m+s}(c_m) \Omega
}\emun
where $U_{st}(\Om\otimes\Om)=\Om$ and $b_1,\dotsc,b_n,c_1,\dotsc,c_m\in\cB,\thinspace n,m\in\N$. Indeed, the mappings $U_{st}$ are unitary. 
The shift is isometric and the unit vector $\Om$ is cyclic which ensures surjectivity. 
Therefore, we may think of the family of Hilbert spaces $(H_t)_{t\geq0}$ as a \hl{tensor product system} in the sense
 of Arveson \cite{Arv}; see Skeide \cite{Ske05b}. In fact, we will see later that is type I.

Let $0\le s = t_0 \leq t_1 \leq \cdots \leq t_{n-1} \leq t_n = t$. 
Using the unitary isomorphism 
$H_{t_0,t_1} \otimes H_{t_1,t_2} \otimes \cdots \otimes H_{t_{n-1},t_n} \cong H_{s,t}$,
in the sequel, we identify
\beq{
\label{TPS}
j_{t_0,t_1}(b_1) \cdots j_{t_{n-1},t_n}(b_n) \Omega = j_{t_0,t_1}(b_1) \Omega \otimes \cdots \otimes j_{t_{n-1},t_n}(b_n) \Omega.
}\eeq

In what follows we will often exploit in an essential way the coalgebra structure of $\ol{\cB}\otimes\cC$ (see Section \ref{prelSEC}) and its interplay with expressions like \eqref{TPS}. The following proposition expresses the core of all such computations. It's proof is an easy verification and we omit it.

\bprop\label{BcCprop}
Let $(\cB,\Delta,\delta)$ and $(\cC,\Lambda,\lambda)$ be coalgebras. Let $D_i$ $(i=1,2)$ be two pre-Hilbert spaces and suppose we have linear mappings $J_i\colon \cB\rightarrow D_i$ and $K_i\colon \cC\rightarrow D_i$. Define the linear functionals $L_i$ on the coalgebra $\ol{\cB}\otimes\cC$ by setting
\beqn{
L_i(\ol{b}\otimes c)
~:=~
\AB{J_i(b),K_i(c)}
}\eeqn
and denote
\baln{
J_1\star J_2
&
~:=~
(J_1\otimes J_2)\circ \Delta
\colon
\cB
~\longrightarrow~
D_1\otimes D_2,
\\
K_1\star K_2
&
~:=~
(K_1\otimes K_2)\circ \Lambda
\colon
\cC
~\longrightarrow~
D_1\otimes D_2.
}\ealn
Then
\beqn{
L_1\star L_2(\ol{b}\otimes c)
~=~
\AB{J_1\star J_2(b),K_1\star K_2(c)}.
}\eeqn
\eprop

Like in Proposition \ref{BcCprop}, in all what follows it is important to pay carefully attention to the several comultiplications of the the coalgebras $\cB$, $\ol{\cB}$, $\cC$, $\ol{\cC}$, $\ol{\cB}\otimes\cC$, and $\ol{\cB}\otimes\cC$, the several convolutions stem from.

\bprop
\label{prop1}
For $c,d\in\cC$ and $T>0$ there exists 
a $C > 0$ such that the following holds.
For each $[s,t]\subset[0,T]$ and $\alpha \in \eZ_{st}$ 
and for each $\beta \in \eZ_{st}$ finer than $\alpha$ we have
\beq{\label{fundest}
\abs{\AB{\vartheta_{\alpha}(c), \vartheta_{\beta}(d)} - e_\star^{(t-s) \psi \circ \vk}(c^* d)} 
< \norm{\alpha}(t-s)C.
}\eeq
\eprop

\lf\lf
\proof
The partitions $\alpha$ and $\beta$ are given by
$\alpha = \{s = s_0  < s_1  < \cdots < s_{l}   = t\}$
and
\baln{
\beta = \{s = s_0  &= t_0^{(1)} < t_1^{(1)} < \cdots < t_{k_1-1}^{(1)}  < t_{k_1}^{(1)}= s_1  \\
& = t_0^{(2)} < t_1^{(2)} < \cdots < t_{k_2-1}^{(2)}  < t_{k_2}^{(2)} = s_2  \\
& \medspace \medspace \thinspace \vdots \\
& = t_0^{(l )} < t_1^{(l )} < \cdots < t_{k_{l }-1}^{(l )} < t_{k_{l }}^{(l )}  = s_{l }    = t\}.  
}\ealn
Denote further $\alpha^{(n)}=\CB{s_{n-1}= t_0^{(n)} < t_1^{(n)} < \cdots < t_{k_{n }-1}^{(n)} < t_{k_{n}}^{(n)}  = s_{n}}$ for $n=1,\ldots, l$. For any pair of partions $\alpha,\beta$ of any interval $[s,t]$ define the linear functionals $L_{\alpha,\beta}$ on $\ol{\cC}\otimes\cC$ by setting $L_{\alpha,\beta}(\ol{c}\otimes d):=\AB{\vt_\alpha(c),\vt_\beta(d)}$. Then, by Proposition \ref{BcCprop},
\beqn{
L_{\alpha,\beta}
~=~
L_{\CB{s_0,s_1},\alpha^{(1)}}\star\ldots\star L_{\CB{s_{l-1},s_l},\alpha^{(l)}}.
}\eeqn
In the concrete form of $L_{\CB{s_{n-1},s_n},\alpha^{(n)}}$ we may rewrite $j_{s_{n-1},s_n}\circ\vk(c)=(j_{t_0^{(n)},t_1^{(n)}}\star\ldots\star j_{t_{k_n-1}^{(n)},t_{k_n}^{(n)}})\circ\vk(c)$, since, by assumption, $(j_{s,t})_{0\leq s\leq t}$ is a Lévy process with respect to the comultiplication of $\cB$. If for any partition $\alpha$ of any interval $[s,t]$ we define the linear functionals $M_\alpha$ on $\ol{\cB}\otimes\cC$ by setting $M_\alpha(\ol{b}\otimes c):=\AB{j_{t_0,t_1}\star\ldots\star j_{t_{n-1},t_n}(b)\Om,\vt_\alpha(c)}$, then, again by Proposition \ref{BcCprop},
\beqn{
L_{\CB{s_{n-1},s_n},\alpha^{(n)}}(\ol{c}\otimes d)
~=~
M_{\CB{t_0^{(n)},t_1^{(n)}}}\star\ldots\star M_{\CB{t_{k_n}^{(n)},t_{k_n-1}^{(n)}}}(\ol{\vk(c)}\otimes d).
}\eeqn
For $\rho\in[0,\norm{\alpha}]$ we define $L_\rho^{(n)}:=L_{\CB{s_{n-1},s_{n-1}+\rho},\alpha^{(n)}(\rho)}$, where
\beqn{
\alpha^{(n)}(\rho)
~:=~
\Bfam{[s_{n-1},s_{n-1}+\rho]\cap\alpha^{(n)}}\cup\bCB{s_{n-1}+\rho}.
}\eeqn
(Roughly speaking, if $\rho\le s_n-s_{n-1}$, then $\alpha^{(n)}(\rho)$ concides with the part of $\alpha^{(n)}$ up to $s_{n-1}+\rho$, and otherwise it adds another interval to the partition.)

We define the linear functionals $M_r:=M_{\CB{\tau,\tau+r}}$ on $\ol{\cB}\otimes\cC$. Note that these do not depend on $\tau\ge0$. We find
\bmun{
M_r(\ol{b}\otimes c)
~=~
M_{\CB{\tau,\tau+r}}(\ol{b}\otimes c)
\\
~=~
\AB{j_{\tau,\tau+r}(b)\Om,j_{\tau,\tau+r}\circ\vk(c)\Om}
~=~
\vp_r(b^*\vk(c))
~=~
\bfam{(\ol{\delta} \otimes \lambda)+ r G+ \eR_r}(\ol{b}\otimes c),
}\emun
where $G(\ol{b}\otimes c):=\psi(b^*\vk(c))$ and $\eR_r$ fulfills the condition of Lemma \ref{lem:lim}. For fixed $[s,t]$, it follows that for every $\ol{c}\otimes d\in\ol{\cC}\otimes\cC$ there exists a constant $C_{c,d}$ such that
\beqn{
\abs{L_\rho^{(n)}(\ol{c}\otimes d)-e_\star^{\rho G}(\ol{\vk(c)}\otimes d)}
~\le~
\snorm{\alpha^{(n)}(\rho)}\rho C_{c,d}
~\le~
\rho^2C_{c,d}
}\eeqn
for all partitions $\alpha^{(n)}$ of $[s_{n-1},s_n]$. (The constant $C_{c,d}$ might depend on $[s,t]$.) From this it is routine to conclude that the $L_\rho^{(n)}$ fulfill the condition of Lemma \ref{lem:lim} at least for all $\ol{c}\otimes d\in\ol{\cC}\otimes\cC$ with the linear first order functional $\ol{c}\otimes d\mapsto\psi\circ\vk(c^*d)$. By takining (finite!) linear combinations, we obtain suitable constants $D_\gamma$ for every $\gamma\in\ol{\cC}\otimes\cC$. From this the statement follows.\qed

\bcor
The net $\bfam{\vt_\alpha(c)}_{\alpha\in\eZ_{st}}$ is a Cauchy net.
\ecor

\proof
We have to show that for $\ve > 0$ there is a $\gamma$ such that
$\alpha , \beta \in \eZ _{st}$,
$\alpha \succ \gamma$ and $\beta \succ\gamma$, implies
$ \norm{\vartheta _{\alpha} \thinspace (c) - \vartheta _{\beta} \thinspace (c)} < \ve $ .
By Proposition \ref{prop1}  there is a $\gamma$ such that for
$\eta \in \eZ _{st}$ with 
$\eta\succ\gamma$, we have
\beq{
\label{approx}
\abs{
\AB {\vartheta _{\gamma} \, (c) , \vartheta _ {\eta} \, (c)} -
e_\star^{(t-s) \psi \circ \vk}(c^* c) } < \frac{\ve ^2}{16} .
}\eeq
So, for $\alpha \succ\gamma$ we have
\beqn{
\begin{split}
\norm{\vartheta _{\eta} \, (c) - \vartheta _{\alpha} \, (c) } ^2 &=  
\AB { \vartheta _{\eta} \, (c) , \vartheta _{\eta } \, (c) } + \AB { \vartheta _{\alpha} \, (c) , \vartheta _{\alpha } \, (c) }  \\
& \qquad - \AB { \vartheta _{\eta} \, (c) , \vartheta _{\alpha } \, (c) }
- \AB { \vartheta _{\alpha} \, (c) , \vartheta _{\eta } \, (c) }  \\
&\leq \frac{\ve ^2}{4} .
\end{split}
}\eeqn
Thus, for $\alpha \succ\gamma$ and $\beta \succ\gamma$ 
\beqn{
\norm{\vartheta _{\alpha} \, (c) - \vartheta _{\beta} \, (c)}
\leq
\norm{\vartheta _{\alpha} \, (c) - \vartheta _{\eta} \, (c)} + \norm{\vartheta _{\beta} \, (c) - \vartheta _{\eta} \, (c)}
\leq \ve. \qedsymbol
}\eeqn
\noqed

\vspace{-2ex}
The limit of the Cauchy net $\bfam{\vt_\alpha(c)}_{\alpha\in\eZ_{st}}$ in $\ol{D_j}$ will
 be denoted by $\vt_{s,t}(c)$.

\brem
Taking the limit of \eqref{fundest} over $\beta\succ\alpha$ for fixed $\alpha$, we find the same estimate for $\AB{\vt_\alpha(c),\vt_{s,t}(c)}$. The fact that \eqref{fundest} does not depend on the precise form of $\alpha$ but only on its width $\norm{\alpha}$ and computations similar to the proof of the corollary, show that $\norm{\vt_\alpha(c)-\vt_{s,t}(c)}$ is small, whenever $\norm{\alpha}$ is sufficiently small. In particular, it follows that 
\beqn{
\lim_{n \to \infty} \vartheta _{\alpha _n} (c) = \vartheta _{s,t}(c)
}\eeqn
for each sequence $\alpha _n$ in $\eZ_{st}$ with 
$\lim _{n \to \infty} \vert \vert \alpha _n \vert \vert = 0$.
\erem

To conclude the proof of Part \ref{trans1} of Theorem \ref{transthm}, we start by observing that 
\beq{\label{vtconv}
\vt_{s,t}(c)
~=~
\vt_{t_0,t_1}\star\ldots\star\vt_{t_{n-1},t_n}.
}\eeq
(To see this, simply take the limit of $\vt_\beta$ over the subnet of partions $\beta\succ\alpha$.) For $\alpha=(s=t_0<t_1<\ldots<t_{n-1}<t_n)\in\eZ_{st}$ ($0\le s<t$) we define
\beqn{
D_{k_\alpha}
~:=~
\ls\BCB{\vt_{t_0,t_1}(c_1)\otimes\ldots\otimes\vt_{t_{n-1},t_n}(c_n )\colon c_1,\ldots,c_n\in\cC}.
}\eeqn
By \eqref{vtconv}, $\vt_{s,t}(c)=\vt_{t_0,t_1}\star\ldots\star\vt_{t_{n-1},t_n}$ it follows $\beta\succ\alpha$ $\Longrightarrow$ $D_{k_\beta}\supset D_{k_\alpha}$. We put $D_{k_{s,t}}:=\bigcup_\alpha D_{k_\alpha}$. Of course, $[s',t']\supset[s,t]$ $\Longrightarrow$ $D_{k_{s',t'}}\supset D_{k_{s,t}}$. We put $D_{k_{t,\infty}}:=\bigcup_{t\le r<s}D_{k_{r,s}}$ and $D_k:=D_{k_{0,\infty}}\ni\Om$. On $D_{k_{s,t}}$ we define an operator by setting
\beqn{
\vartheta_{t_0,t_1}(c_1) \otimes \cdots \otimes \vartheta_{t_{n-1},t_n}(c_n)
~\longmapsto~
\vartheta_{t_0,t_1}(c_{(1)}c_1) \otimes \cdots \otimes \vartheta_{t_{n-1},t_n}(c_{(n)}c_n).
}\eeqn
To see that this is well-defined, we simply observe that the the operator has a formal adjoint on that domain, namely, simply the operator whit $c$ replaced by $c^*$. (By taking joint refinements, if necessary, we may assume that the two vectors we choose to check the adjoint condition are in the same $D_{k_\alpha}$.) We extend this operator by amplification to an operator $k_{s,t}(c)$ on $D_k=D_{k_{0,s}}\otimes D_{k_{s,t}}\otimes D_{k_{t,\infty}}$. Clearly, $c\mapsto k_{s,t}(c)$ is multiplicative, so that the $k_{s,t}$ define a family of \nbd{*}homomorphisms. A simple application of coassociativity (and, once more, \eqref{vtconv}) shows that $k_{r,s}\star k_{s,t}=k_{r,t}$ for $r<s<t$. Therefore,  the family of mappings $k_{s,t}$ forms a Lèvy process on $\cC$ over $(D_k,\Om)$ with generator $\psi\circ\vk$. That $D_k$ is dense in $H_k$, will follow from the proof of Part \ref{trans2}.

\subsection{Proof of Part \ref{trans2} of Theorem \ref{transthm}}\label{p2SS}

By Part \ref{trans1} of Theorem \ref{transthm} we know that the $\zeta_\alpha$ converge in norm to something that determines a Lévy process $\tilde{j}$ on $D_{\tilde{j}}$ that is equivalent to $j$. In particular, $\AB{\zeta_\alpha(b),\zeta_\alpha(b)}~\to~e_\star^{(t-s)\psi}(b^*b)$ $=$ $\AB{j_{s,t}(b)\Om,j_{s,t}(b)\Om}$. Therefore, the only thing that remains to be shown in order to see that $\snorm{\zeta_\alpha-j_{s,t}(b)\Om}^2$ $\to$ $0$, in other words, that $\tilde{j}=j$, is the following proposition.

\bprop
\label{prop2}
For all $b,d\in\cB$ we have
\[
\lim_\alpha\AB{\zeta_\alpha(b), j_{s,t}(d)\Om}
~=~
e_\star^{(t-s) \psi}(b^* d).
\]
\eprop

\proof
Let $\alpha=\CB{s=t_0<t_1<\ldots<t_{n-1}<t_n=t}$ and write $j_{s,t}=j_{t_0,t_1}\star\ldots\star j_{t_{n-1},t_n}$. Then, as in the proof of Proposition \ref{prop1}, from Proposition \ref{BcCprop} we find
\beqn{
\AB{\zeta_\alpha(b), j_{s,t}(d)\Om}
~=~
L_{t_1-t_0}\star\ldots\star L_{t_n-t_{n-1}}(\ol{b}\otimes d),
}\eeqn
Where we define the linear functionals $L_r(\ol{b}\otimes d):=\AB{k_{0,r}\circ\wt{\vk}(b)\Om,j_{0,r}(d)\Om}$ on $\ol{\cB}\otimes\cB$.

We are done, if we show the the $L_r$ fulfill the conditions of Lemma \ref{lem:lim} with the correct linear term. In fact, if in \eqref{fundest} we insert $\alpha=\CB{0,r}$ (so that $\norm{\alpha}=r$) and perform the limit over $\beta$, the estimate remains valid for $\AB{\vt_{\CB{0,r}}\circ\wt{\vk}(d),k_{0,r}\circ\wt{\vk}(b)\Om}=\ol{L_r(\ol{b}\otimes d)}$.\qed

\lf
This ends also the proof of Part \ref{trans2} of Theorem \ref{transthm}.

\bcor
The vectors $k_{s,t}\Om,\thinspace c\in\cC,$ generate $\ol{D_j}$ in the sense that
\bmun{
\ol{D_j}=\ol{D_k}=\cls\BCB{k_{t_0,t_1}(c_1) \cdots k_{t_{n-1},t_n}(c_n) \Omega \colon\\
n \in \mathbb{N}, \thinspace 0 \leq s \leq t < \infty, \thinspace s = t_0 \leq t_1 \leq \cdots \leq t_{n-1} \leq t_n=t, \thinspace c_1, \dotsc ,c_n \in \mathcal{C}}.
}\emun
\ecor

\section{Applications of the transformation theorem}\label{applSEC}

\subsection{Realization of quantum Lévy processes on Boson Fock space}\label{LPFockSEC}

Now we apply the Transformation Theorem (Theorem \ref{transthm}) to the Example \ref{exgl}. 
To that goal, let $(\cB,\Delta,\delta)$ be some \nbd{*}bialgebra and let $\bfam{j_{s,t}}_{0\leq s\leq t<\infty}$ be 
a cyclic 
Lévy process on $\cB$ over $(D_j,\Om)$ with generator $\psi$. 
In view of Part \ref{trans1} of Theorem \ref{transthm} we have that
\beqn{
k_{s,t}(\hat b) \Om := \lim j_{t_0,t_1}(b)\dotsc j_{t_{n-1},t_n}(b)\Om
}\eeqn
for $b\in\cB_1$ defines a cyclic Lévy process $\bfam{k_{s,t}}_{0\leq s\leq t<\infty}$ on
 $(\C(\cB_1),\Lambda,\lambda)$ over $(D_k,\Om)$ where $D_k$ is a linear subspace of $\ol{D_j}$. 
Thus, for each pair $k_{s,t}(\hat b)\Om,\thinspace k_{s,t}(\hat c)\Om$ for $b,c\in\cB_1$ and $0\leq s\leq t<\infty$
 we have
\beqn{
\AB{k_{s,t}(\hat b) \Omega, k_{s,t}(\hat c) \Omega} = e^{(t-s)\psi(b^*c)}.
}\eeqn
The generator $\psi$ defines a coboundary by \eqref{korand}. 
Thus, we compute
\beqn{
\begin{split}
\AB{e^{-(t-s)\psi(b)}k_{s,t}(\hat b)\Om,e^{-(t-s)\psi(c)}k_{s,t}(\hat c)\Om}
&=e^{(t-s)(-\psi(b^*)-\psi(c)+\psi(b^*c))}\\
&=e^{(t-s)\AB{\eta(b),\eta(c)}}\\
&=\AB{E(\eta(b)\otimes\eins _{[s,t]}),E(\eta(c)\otimes\eins _{[s,t]})}
\end{split}
}\eeqn
where $\eta\colon\cB_1\rightarrow K$ is the canonical mapping to 
a dense linear subspace $K$ of a Hilbert space $\ol{K}$ and $E(\eta(\cdot)\otimes\eins _{[s,t]})$ denotes
 the exponential vector of $\eta(\cdot)\otimes\eins _{[s,t]}$ in the Boson Fock space $\Gamma_s(\RL^2([s,t], \ol{K}))$.
 Here $\eta(\cdot)\otimes\eins _{[s,t]}$ denotes the function in $\RL^2([s,t], \ol{K})$ which is a constant
equal to $\eta(\cdot)$
 on the interval $[s,t]$ and zero elsewhere. The space $K$ is obtained by applying a GNS-type construction to $\psi$. 
Hence, 
\beqn{
k_{s,t}(b)\Om\cong e^{(t-s)\psi(b)}E(\eta(b)\otimes\eins _{[s,t]})\in\Gamma_s(\RL^2([s,t],\ol{K}))
}\eeqn
where $b\in\cB_1, \psi(b)\in\C$ and $\eta(b)\in K$.
In other words, the vectors $k_{s,t}(b)\Om$ behave like exponential vectors 
in the Boson Fock space $\Gamma_s(\RL^2([s,t],\ol{K}))$. 
Moreover, the vectors $k_{s,t}(b)\Om$ \lq generate\rq \
 the Hilbert subspace $D_{k_{s,t}}$ of $\ol{D_k}$ where 
\[D_{k_{s,t}}=\ls\BCB{k_{t_0,t_1}(c_1) \cdots k_{t_{n-1},t_n}(c_n) \Omega \colon
n \in \mathbb{N},\thinspace s = t_0 \leq t_1 \leq \cdots \leq t_{n-1} \leq t_n=t, \thinspace c_1, \dotsc ,c_n \in \mathcal{C}}.\]
Therefore, we have $\ol{D_{k_{s,t}}}\cong\Gamma_s(\RL^2([s,t],\ol{K}))$ and thus $\ol{D_k}\cong\Gamma_s(\RL^2(\R^+,\ol{K}))$.
Part \ref{trans2} of the Transformation Theorem states that the vectors $k_{s,t}(b)\Om,\thinspace b\in\cB_1,$ are total in
 $\ol{D_{j_{s,t}}}\subset\ol{D_j}$ as well, $\rm{i.e.}$,
\beqn{
\ol{D_j}=\ol{D_k}\cong\Gamma_s(\RL^2(\R^+,\ol{K})).
}\eeqn
So we proved that each cyclic quantum Lévy process on a \nbd{*}bialgebra can
 be realized on a Boson Fock space $\Gamma_s(\RL^2(\R^+,\ol{K}))$.

\subsection{Construction of quantum Lévy processes}

In the situation of Section \ref{LPFockSEC}, an application of Part \ref{trans2} of Theorem \ref{transthm} allows to reconstruct $j_{s,t}$ from the process $k_{s,t}$ on the group-like \nbd{*}bi\-al\-ge\-bra. The realization of the latter on the Fock space can simply be written down. In the present section we describe a realization on the Fock space that rather parellels the construction in \cite{MSchue93} with the help of quantum stochastic calculus.

We will describe the construction of $k_{s,t}$ out of $j_{s,t}$ in part
\ref{trans1} of the transformation theorem
by the short hand writing
\beq{\label{TF}
{\smash{\prod^{\rightarrow}}\rule[-1.5ex]{0pt}{4.3ex}}_{\GL}^\star (j_{s,t} \circ \vk ) = k_{s,t} .
}\eeq
We call $k_{s,t}$ the \hl{infinitesimal convolution product} of $j_{s,t} \circ \vk$.

Applying our result to the situation of Example \ref{exta} and \ref{extarev} with $\vk = \id$ there are two possibilities.
If we put $\SB$ equal to the induced 
\nbd{*}bialgebra and $\SC$ equal to the primitive 
\nbd{*}bialgebra, then
for $b \in \SB _0$ we have
\beq{
\label{add}
\vartheta _{\ga } (b) = 
 \sum_{i = 1} ^n j_{t_{i-1} , t_i} (b)\Om
}\eeq
and Part \ref{trans1} of Theorem \ref{transthm} tells us that \ref{add} converges to
\beqn{
I_{s,t} (b) = A_{s,t} (\eta (b ^* )) + \GL _{s,t} ( \rho (b)) +
A^* _{s,t} (\eta (b)) + \psi (b) \, (t - s)
}\eeqn
in norm where 
$A_{s,t} , \GL _{s,t} , A ^ * _{s,t}$
denote the annihilation, preservation and creation operators
of the interval $[s, t]$ on Boson Fock space
$\GG _s (\RL ^2 (\R _+ , \ol{ K }))$; see the preceding section. For arbitrary $b\in\cB$ we find
\beqn{
I_{s,t} (b) =\delta(b)I+ A_{s,t} (\eta (b ^* )) + \GL _{s,t} ( \rho (b)-\delta(b)) +
A^* _{s,t} (\eta (b)) + \psi (b-\delta(b)) \, (t - s).
}\eeqn
$I_{s,t}$ is the \hl{additive generator process} of the Lévy process $j_{s,t}$. (It is addiditve on $\cB_0$, repcetively, the process $I_{s,t}-\delta I$ is additive.)

We construct $j_{s,t}$ out of $I_{s,t}$ if we take the primitive 
\nbd{*}bialgebra
for $\SB$ and the induced one for $\SC$.
Then by Part \ref{trans1} of Theorem \ref{transthm} we obtain
$j_{s,t}$ as the limit
\beqn{
j_{s,t} = {\smash{\prod^{\rightarrow}}\rule[-1.5ex]{0pt}{4.2ex}}_{\ST (\GD _0)}^\star I_{s,t} 
}\eeqn
of the convolution products of the generator process
where now, of course, convolution is with respect to the original comultiplication
$\GD$ of $\SB$.
So our procedure allows, like quantum stochastic calculus, a construction of the
Lévy process $j_{s,t}$ from the elementary processes
$A_{s,t}, \GL _{s,t}, A^* _{s,t}$ on Boson Fock space. In fact, if $dt$ is ``small'', then in all relevant formula one may substitute $j_{t,t+dt}$ with $I_{t,t+dt}$. We find
\beqn{
j_{s,t+dt}-j_{s,t}
~=~
j_{s,t}\star j_{,t+dt}-j_{s,t}
~\simeq~
j_{s,t}\star I_{t,t+dt}-j_{s,t}
~=~
j_{s,t}\star(I_{t,t+dt}-\delta I).
}\eeqn
If we put $dI_t=I_{s,t+dt}-I_{s,t}$ (independent of $s<t$), this gives an immediate meaning to 
\beqn{
j_{s,t}
~=~
\delta I+\int_s^t j_{s,r}\star dI_r
}\eeqn
as a quantum stochastic integral. We remark that this interpretation as integral is not limited to the above choice. Whenever $k$ is a transformed process obtained from $j$ via \eqref{TF}, then it fulfills
\beqn{
k_{s,t}
~=~
\delta I+\int_s^t k_{s,r}\star(dj_t\circ\vk),
}\eeqn
where $dj_t:=j_{t,t+dt}-\delta I$.

\subsection{Classical Lévy processes and unitary evolutions}

Let $G$ be a topological group and denote by $\SR (G)$ the space of all coefficient
functions of continuous finite-dimensional representations of $G$.
Then $f \in \SR (G)$ iff there are $n \in \N$ and continuous complex-valued functions $f_1 , \dots f_n , g_1 , \dots g_n$ on $G$
such that
\beqn{
f(x y) = \sum_{i = 1} ^n f_i (x) \, g_i (y)  \  \forall x, y \in G .
}\eeqn
$\SR (G)$ is a commutative 
\nbd{*}algebra.
By setting
\beqn{
\GD f = \sum_{i = 1 } ^n f_i \ot g_i  , \  \
\gd f = f(e) 
}\eeqn
$\SR (G)$ becomes a commutative Hopf 
\nbd{*}algebra.
In various cases (e.g., when $G$ is compact or locally compact abelian) 
the group $G$ is uniquely determined by $\SR (G)$.
Let us assume that $G$ is compact.
Then $\SR (G)$ is the Kre\^in algebra of $G$. A classical Lévy process
$X_t$ on $G$ gives rise to a quantum Lévy process $j_t$ on $\SR (G)$ 
by putting
$j_t (f) = f \circ X_t$.
Here
$j_t = j_{0t}$ and $j_{s,t} = (j_s \circ S ) \cst j_t$ where $S$ is the antipode of $\SR (G)$.
Let us specialize to the case when $G$ is the group $\SU _d$ of unitary $d \times d$-matrices.
Then $\SR (G)$ equals the Hopf 
\nbd{*}algebra $\C [ x_{kl}, x_{kl} ^ *  ; k, l = 1, \dots d ]$ divided by the 
\nbd{*}ideal generated by the elements which are the entries of the matrices 
$x \, x ^ *  - \eins $ and $x ^ *  \, x - \eins$ where we put $x = (x_{kl} )_{k, l = 1, \dots d}$.
The comultiplication is given by
$\GD x_{kl} = \sum_{i = 1} ^d x_{ki} \ot x_{il}$
and the counit by $\gd x_{kl} = \gd _{kl}$.
The antipode is given by $S (x_{kl}) = x_{lk} ^ * $.
By replacing the commuting indeterminates $x_{kl}$ by non-commuting indeterminates,
we define a non-commutative 
\nbd{*}bialgebra
\beqn{
\C \langle x_{kl}, x_{kl} ^ *  ; k, l = 1, \dots d  \rangle / x \, x^ *  = \eins, x^ *  \, x = \eins 
}\eeqn
which we denote by $\SU \langle d \rangle$. (It is easy to see that $\SU \langle d \rangle$ is not a Hopf algebra.)
Lévy triples on $\SU \langle d \rangle$ are given by a Hilbert space $\ol{K}$, a unitary operator $W$ on
$\C ^d \ot \ol{K}$, a matrix $L \in M_d (\C ) \ot \ol{K}$ and a self-adjoint matrix $H \in M_d (\C )$
via the equations
\beqn{
\begin{split}
\rho (x_{kl} ) &= W_{kl} \in \SB (\ol{K} )   \\
\eta (x_{kl} ) &= L_{kl}   \\
\psi (x_{kl} ) &= - \frac12 (L L^ *  )_{kl} + \ri \, H_{kl} ;
\end{split}
}\eeqn
cf. \cite{MSchue93}.
The generator process is given by matrices $ \SI _{s,t} \in M_d (\C ) \ot 
\GG (\RL ^2 (\R _+  , \ol {K} ))$ with
\beqn{
(\SI _{s,t} )_{ij} = - A_{s,t}  ( (W ^ *  L) _{ji} ) +
\GL _{s,t} ((W - \eins ) _{ij} ) + A_{s,t}^* (L_{ij} ) 
+ (\ri \, H - \frac{1}{2} (L L^ *  ) )_{ij} \, (t - s)
}\eeqn
The transformation Theorem \ref{transthm} says that
\beqn{
\SI _{t_0 , t_1 } \SI _{t_1 , t_2 }  \dots \SI _{t_{n - 1} , t_n}
}\eeqn
converges to the Lévy process $U_{s,t}$ which  is the unitary process
on $\C ^d \ot \GG (\RL ^2 (\R _+ , \ol{k} ))$ given by
$(U_{s,t} )_{ij}  = j_{s,t} (x_{ij})$.
This is a generalization of a construction already given in \cite{vW84}.
A classical Lévy process on $\SU _d$ is a special case of a QLP on $\SU \langle d \rangle$.

\subsection{Azéma martingales}\label{AzemaSEC}

Consider the 
\nbd{*}algebra $\C \langle x, x ^* , y \rangle$ generated by $x$ and a self-adjoint $y$.
For $q \in \R$ divide $\C \langle x, x ^* , y \rangle$ by the 
\nbd{*}ideal generated by the element
$x y - q y x$
to obtain a 
\nbd{*}algebra $\SA$.
On $\SA$ we consider two 
\nbd{*}bialgebra structures.
The first is the one with $x$ (and $x ^*$) primitive and with $y$ group-like,
the second is given by
\beqn{
\begin{split}
\GD x &= x \ot y + \eins \ot x \  \text{and}  \  \ \gd x = 0   \\
\GD y &= y \ot y \  \text{and}  \ \  \gd y = 1
\end{split}
}\eeqn
and maybe called the \hl{Azéma
\nbd{*}bialgebra for parameter} $q$.
Again we apply our results to these two 
\nbd{*}bialgebras with $\vk = \id$.
If we choose for generator
\beqn{
\psi (M(x, x^* ) \, y ^k ) = \left\{
\begin{array}{ll}
1 & \mbox{if} \ M(x, x^* ) = x x ^*  \\
0 & \mbox{otherwise}
\end{array}
\right.
}\eeqn
$M(x, x ^* ) \in \SA$ a monomial in $x$ and $x^*$, $ k \in \N _0$, then
$K = \C$, $\eta (x ^ * ) = 1, \ \eta (x) = 0$, $\rho (x) = 0$ and $\rho (y) = q$.
The linear functional $\psi$ is the generator of the quantum $q$-Azéma martingale $(X_t , X_t ^* , Y_t )$
if we consider the Azéma
\nbd{*}bialgebra, and it generates the process
$(A_t , A^* _t , Y_t )$ in the case of the primitive/group-like structure of $\SA$
where
$Y_t$ is the second quantization of multiplication by $q \, \eins _{[0, t ]}$.
The process $X_t$ satisfies the quantum stochastic differential equation
\beqn{
\rd X_t = (q - 1) X_t \, \rd \GL _t + \rd A_t , \ X_0 = 0 ; 
}\eeqn
see \cite{Partha, MSchue93}
An application of Part \ref{trans1} of Theorem \ref{transthm}
yields the formulae
\beqn{
W _t = \lim \sum_{j = 0}^{n - 1} Z_{t_j , t_{j + 1}}
}\eeqn
and
\beqn{
Z_t = \lim
\bigl( W_{t_0 , t_1}\, Y_{t_1 , t_2} \dots
Y_{t_{n - 1}, t_n}  +
W_{t_1 t_2 } \, Y_{t_2 , t_3} \dots Y_{t_{n - 1} , t_n} + \dots 
+ W_{t_{n - 2} , t_{n - 1}} \, Y_{t_{n - 1}, t_n} + W_{t_{n - 1} , t_n}
\bigr)
}\eeqn
where $W_t$ and $Z_t$ denote the Wiener process and the $q$-Azéma martingale on Boson Fock space respectively.

\setlength{\baselineskip}{2.5ex}


\begin{thebibliography}{xx}

\vspace{-.5ex}
\bibitem{ASchvW} L. Accardi, M. Schürmann, W. von Waldenfels: Quantum independent increment 
processes on superalgebras. Math. Z. 
\bf {198}, 451-477 (1988)

\bibitem{Arv}
W. Arveson: Noncommutative Dynamics and A-Semigroups. Springer Monographs in Matmematics. 
New York Berlin Heidelberg, Springer 2003

\bibitem{Fra06} U. Franz: Lévy processes on quantum groups and dual groups.
In: Franz, U, Schürmann, M. (eds.) Quantum independent increment processes II.
Lect. Notes Math., vol. 1866. New York Berlin Heidelberg, Springer 2006

\bibitem{HuPa84} R.L. Hudson, K.R. Parthasarathy: Quantum Ito's formula 
and stochastic evolutions. Commun. Math. Phys. \bf {93}, 301-323
(1984)

\bibitem{LiSk05}
V. Liebscher and M. Skeide: Constructing units in product systems. Proc.\ Amer.\ Math.\ Soc. \bf{136}, 989--997, (2008),
electronically Nov 2007, (ar\-Xiv: math.OA/0510677).

\bibitem{Partha} K.R. Parthasarathy: Azéma martingales and quantum stochastic calculus. In: Bahadur, R.R. (ed.)
Proc. R.C. Bose Memeorial Symposium.
New Delhi, Wiley Eastern 1990

\bibitem{PaSu98}
K.R. Parthasarathy and V.S. Sunder:
Exponentials of indicater functions are total in the boson Fock
  space $\G(L^2[0,1])$.
In: R.L. Hudson and J.M. Lindsay, editors, Quantum Probability
  Communications X, number~X in QP-PQ:, pages 281--284. World Scientific,
  1998.

\bibitem{MSchue93} M. Schürmann: White Noise on Bialgebras.
 Springer Lect.  Notes Math., vol. 1544. New York Berlin Heidelberg, Springer 1993

\bibitem{Ske00a}
M. Skeide:
Indicator functions of intervals are totalizing in the symmetric
  Fock space $\Gamma(L^2(\R_+))$.
In L. Accardi, H.-H. Kuo, N. Obata, K. Saito, {Si Si}, and L. Streit,
  editors, Trends in contemporary infinite dimensional analysis and
  quantum probability, volume~3 of Natural and Mathematical Sciences
  Series, pages 421--424. Istituto Italiano di Cultura (ISEAS), Kyoto, 2000.
Volume in honour of Takeyuki Hida, (Rome, Volterra-Pre\-print
  1999/0395).

\bibitem{Ske05b}
M. Skeide: ~
Lévy processes and tensor product systems of Hilbert modules.
In M.~Sch\"urmann and U.~Franz, editors, {\em Quantum Probability and
Infinite Dimensional Analysis --- From Foundations to Applications}, number
XVIII in Quantum Probability and White Noise Analysis, pages 492--503. World
Scientific, 2005.

\bibitem{vW84} W. von Waldenfels: Ito solution of the linear quantum stochastic 
differential equation describing light emission
and absorption. In: Accardi, L., Frigerio, A., Gorini, V. (eds.) Quantum probability and applications to the theory of
irreversibel processes. Proceedings, Villa Mondragone 1982.
Lect. Notes Math., vol. 1055. New York Berlin Heidelberg, Springer 1984



\end{thebibliography}


\noindent
Michael Schürmann: {\small\itshape Institut für Mathematik und Informatik}, {\small\itshape 
Ernst-Moritz-Arndt-Universität Greifswald}, {\small\itshape 17487 Greifswald, Germany},
{\small\itshape E-mail: \tt{schurman@uni-greifswald.de}},
\\
{\small{\itshape Homepage: \tt{http://www.math-inf.uni-greifswald.de/algebra}}}

\vspace{-1ex}\lf\noindent
Michael Skeide: {\small\itshape Dipartimento S.E.G.e S.}, {\small\itshape Università degli Studi del Molise}, {\small\itshape Via de Sanctis}, {\small\itshape 86100 Campobasso, Italy}, {\small{\itshape E-mail: \tt{skeide@math.tu-cottbus.de}}},
\\
{\small{\itshape Homepage: \tt{http://www.math.tu-cottbus.de/INSTITUT/lswas/\_skeide.html}}}

\vspace{-1ex}\lf\noindent
Silvia Volkwardt: {\small\itshape Institut für Mathematik und Informatik}, {\small\itshape 
Ernst-Moritz-Arndt-Universität Greifswald}, {\small\itshape 17487 Greifswald, Germany},
{\small\itshape E-mail: \tt{svolkwardt@hotmail.com}}


\end{document}